\documentclass[11pt]{amsart}
\usepackage{amscd}      
\usepackage{amsbsy, amsthm, eucal, graphics,epic,amsbsy, stmaryrd, epic, eepic}
\usepackage{amssymb}
\usepackage{amsmath}
\usepackage[all]{xy}
\usepackage[dvips]{color}

\usepackage{latexsym}
\usepackage{epsfig}
\usepackage{amsfonts}
\usepackage{enumerate}

\newtheorem{prop}{Proposition}[section]
\newtheorem{lem}[prop]{Lemma}

\newtheorem{cor}[prop]{Corollary}
\newtheorem{thm}[prop]{Theorem}
\newtheorem{rmk}[prop]{Remark}

\newtheorem{notation}[prop]{Notation}

\newtheorem{defn}[prop]{Definition}

\newtheorem{property}[prop]{Property}
\newtheorem{example}[prop]{Example}

\newenvironment{pf}{\begin{trivlist}\item[]{\sc Proof.}}%
            {\nolinebreak $\Box$ \end{trivlist}}

\newcommand{\noprint}[1]{}

\renewcommand{\tilde}{\widetilde}

\newcommand{\YY}{{\mathfrak Y}}

\newcommand{\MM}{{\mathfrak M}}

\newcommand{\ldiag}[1]%
       {\makebox[0cm]{${\scriptstyle#1}\downarrow\phantom{\scriptstyle#1}$}}
\newcommand{\ldiagup}[1]%
       {\makebox[0cm]{${\scriptstyle#1}\uparrow\phantom{\scriptstyle#1}$}}
\newcommand{\rdiag}[1]%
       {\makebox[0cm]{$\phantom{\scriptstyle#1}\downarrow{\scriptstyle#1}$}}
\newcommand{\sediagr}[1]%
       {\makebox[0cm]{$\phantom{\scriptstyle#1}\searrow{\scriptstyle#1}$}}
\newcommand{\nediagr}[1]%
       {\makebox[0cm]{$\phantom{\scriptstyle#1}\nearrow{\scriptstyle#1}$}}
\newcommand{\rdiagup}[1]%
       {\makebox[0cm]{$\phantom{\scriptstyle#1}\uparrow{\scriptstyle#1}$}}
\newcommand{\swdiag}[1]%
       {\makebox[0cm]{$\phantom{\scriptstyle#1}\swarrow{\scriptstyle#1}$}}
\newcommand{\sediag}[1]%
       {\makebox[0cm]{${\scriptstyle#1}\searrow\phantom{\scriptstyle#1}$}}
\newcommand{\nediag}[1]%
       {\makebox[0cm]{${\scriptstyle#1}\nearrow\phantom{\scriptstyle#1}$}}

\newcommand{\doublearrowstack}[2]%
                      {{{{\scriptstyle#1}\atop{\textstyle\longrightarrow}}\atop{{\textstyle\longrightarrow}\atop{\scriptstyle#2}}}}
\newcommand{\rightleftarrowstack}[2]%
                      {{{{\scriptstyle#1}\atop{\textstyle\longrightarrow}}\atop{{\textstyle\longleftarrow}\atop{\scriptstyle#2}}}}
\newcommand{\leftrightarrowstack}[2]%
                      {{{{\scriptstyle#1}\atop{\textstyle\longleftarrow}}\atop{{\textstyle\longrightarrow}\atop{\scriptstyle#2}}}}

\newcommand{\overtoparrow}%
{\makebox[0cm]{\beginpicture \setcoordinatesystem units
<.8cm,.4cm> point at 0 0 \setplotarea x from -3 to 3, y from 0 to
1 \setquadratic \plot -3 0 0 1 3 0 / \put{\vector(3,-1){0}}[Bl] at
3 0
\endpicture}}

\newcommand{\underbottomarrow}%
{\makebox[0cm]{\beginpicture \setcoordinatesystem units
<.8cm,.4cm> point at 0 0 \setplotarea x from -3 to 3, y from 0 to
1 \setquadratic \plot -3 1 0 0 3 1 / \put{\vector(3,1){0}}[Bl] at
3 1
\endpicture}}

\newcommand{\ses}[5]%
{0\longrightarrow#1\stackrel{#2}{ \longrightarrow}#3\stackrel{#4}{
\longrightarrow}#5\longrightarrow0}

\newcommand{\dt}[6]%
{#1\stackrel{#2}{longrightarrow}#3
\stackrel{#4}{\longrightarrow}#5 \stackrel{#6}{\longrightarrow}
#1[1]}

\newcommand{\cat}[1]%
{(\mbox{\rm #1})}

  \def\bbC{{\mathbb C}}  

\def\bbG{{\mathbb G}}    
\def\bbQ{{\mathbb Q}} \def\bbZ{{\mathbb Z}}   \def\bbN{{\mathbb N}}

\newcommand{\clD}{{\mathcal{D}}}
\newcommand{\clC}{{\mathcal{C}}}

\newcommand{\clG}{{\mathcal{G}}}
\newcommand{\clF}{{\mathcal F}}
\newcommand{\clK}{{\mathcal{K}}}
\newcommand{\clL}{{\mathcal{L}}}
\newcommand{\clM}{\mathcal{M}}
\newcommand{\clN}{\mathcal{N}}

\newcommand{\clS}{\mathcal{S}}

\newcommand{\clR}{{\mathcal R}}
\newcommand{\clT}{{\mathcal{T}}}

\newcommand{\clO}{{\mathcal{O}}}
\newcommand{\clQ}{{\mathcal{Q}}}
\newcommand{\clX}{{\mathcal{X}}}
\newcommand{\clY}{{\mathcal{Y}}}

\def\longto{\longrightarrow}
\def\isomto{\stackrel{\sim}{\longto}}
\providecommand{\abs}[1]{\lvert#1\rvert}
\newcommand{\Spec}{\operatorname{Spec}}
\newcommand{\Pic}{\operatorname{Pic}}
\newcommand{\thickslash}{\mathbin{\!\!\pmb{\fatslash}}}

\newcommand{\barM}{{\overline{M}}}
\newcommand{\ttau}{{\tilde{\tau}}}

\begin{document}
\title[GW theory of banded gerbes over schemes]{Gromov-Witten theory\\ of banded gerbes over schemes}
\author{Elena Andreini}
\address{International School for Advanced Studies\\Via Beirut 2-4\\
34100 Trieste\\ Italy}
\email{andreini.elena@gmail.com}
\author{Yunfeng Jiang}
\address{Department of Mathematics\\ Imperial College London\\ South Kensington Campus\\ London SW7 2AZ\\ United Kingdom}
\email{y.jiang@imperial.ac.uk}
\author{Hsian-Hua Tseng}
\address{Department of Mathematics\\ Ohio State University\\ 100 Math Tower, 231 West 18th Ave.\\Columbus\\ OH 43210\\ USA}
\email{hhtseng@math.ohio-state.edu}

\date{\today}
\begin{abstract}
Let $X$ be a smooth complex projective algebraic variety. Let $\clG$ be a $G$-banded gerbe with $G$ a finite abelian group.  We prove an exact formula expressing genus $g$ orbifold Gromov-Witten invariants of $\clG$ in terms of those of $X$. 
\end{abstract}


\maketitle

\tableofcontents
\section{Introduction}
This paper is a sequel to our paper \cite{AJT09-gerbes}. Our goal is to study the Gromov-Witten theory of \'etale abelian banded gerbes over varieties. Let $X$ be a smooth projective variety over $\mathbb{C}$ and let $G$  be a finite abelian group (viewed as a constant group scheme over $X$). Let $\epsilon: \clG\to X$ be a $G$-banded gerbe over $X$.  As explained in \cite{AJT09-gerbes} our study of Gromov-Witten theory of $\clG$ is motivated by the so-called {\em decomposition conjecture} in physics \cite{PaSha06Clust}. In the case of a banded  gerbe $\clG \to X$ the decomposition conjecture states that the Gromov-Witten theory of $\clG$ is equivalent to the Gromov-Witten theory of the disjoint union of $\abs{\hat{G}}$ copies of $X$, where $\abs{\hat{G}}$ denotes the number of irreducible representations of $G$. See \cite{TT} for more discussions on the mathematical aspects of the decomposition conjecture.

In \cite{AJT09-gerbes} we proved this conjecture for genus $0$ Gromov-Witten theory of $\clG$. Our approach consists of two steps. First, we studied the natural morphism\footnote{This natural morphism is obtained by composing stable maps to $\clG$ with the map $\omega: \clG\to X$ and taking the relative coarse moduli space morphism.} $p: \clK_{0,n}(\clG,\beta)\to\overline{M}_{0,n}(X, \beta)$ between the moduli spaces of stable maps to $\clG$ and to $X$. We showed that this morphism can be factored as $\clK_{0,n}(\clG,\beta)\to P_n\to \overline{M}_{0,n}(X, \beta)$, where $\clK_{0,n}(\clG,\beta)\to P_n$ is a $G$-gerbe and $P_n\to \overline{M}_{0,n}(X, \beta)$ is obtained as a base-change of  a functorial construction \`a la Matsuki-Olsson \cite{MaO05}  applied to the stack  $\mathfrak{M}_{0,n}$ of genus $0$ prestable pointed curves\footnote{Given an algebraic stack $X$ endowed with a locally free log structure the construction associates to it an algebraic stack $\clX$ defined as the category fibered in groupoids whose objects are morphisms to $f:T\to X$ and simple morphisms of locally free log structures $f^*\clM\to \clM'$.}. Such a construction produces a stack over  $\mathfrak{M}_{0,n}$ which can be identified with an open substack of the stack of twisted curves $\MM_{0,n}^{tw}$. This result allows us  to compare the push-forward $p_*[\clK_{0,n}(\clG,\beta)]^{vir}$ of the virtual fundamental class with the virtual fundamental class $[\overline{M}_{0,n}(X,\beta)]^{vir}$. This in turn gives a comparison of genus $0$ Gromov-Witten invariants of $\clG$ and $X$. In the second step we use some standard finite group theory to interpret our comparison result as a statement in agreement with the decomposition conjecture.

In this paper we study the {\em higher genus} Gromov-Witten invariants of $\clG$ following a similar strategy. Again we begin with studying the natural morphism
$$p: \clK_{g,n}(\clG,\beta)\to\overline{M}_{g,n}(X, \beta)$$
between moduli spaces of stable maps. We first show that $p$ can be similarly factored as $$\clK_{g,n}(\clG,\beta)\to P_{g,n}\to \overline{M}_{g,n}(X, \beta),$$ where $P_{g,n}$ is an open substack of $\MM_{g,n}^{tw}\times_{\MM_{g,n}} \overline{M}_{g,n}(X, \beta)$.
Even though, unlike in the genus $0$ case,    $\clK_{g,n}(\clG,\beta)\to P_{g,n}$ is not as simple as the structure morphism of a gerbe,  we show that  it is \'etale (see Proposition \ref{t-etaleness-prop}). It turns out that this property together   with the properties of $P_{g,n}\to \overline{M}_{g,n}(X, \beta)$ are sufficient to prove a comparison result between classes $p_*[\clK_{g,n}(\clG,\beta)]^{vir}$  and $[\overline{M}_{g,n}(X,\beta)]^{vir}$, see Theorem \ref{push-thm}. Using this we then draw conclusions on higher genus Gromov-Witten invariants by the method identical to the genus $0$ case, see Theorem \ref{decomp_any_genus}.

\subsection*{Conventions}
Unless otherwise mentioned, we work over $\bbC$ throughout this paper. By an {\em algebraic stack}  we mean an  algebraic stack over $\bbC$ in the sense of \cite{Art74}.  By a {\em Deligne-Mumford stack} we mean an algebraic stack over $\bbC$ in the sense of \cite{DM69}. We assume moreover that all stacks (and schemes) are locally noetherian and locally of finite type. Following \cite{KatoLog}, logarithmic structures  are considered on the \'etale site of schemes. For the  extension of logarithmic structures to stacks, see \cite{OlssLog03}. Given a scheme (or a stack) $X$, a geometric point $x$ of $X$, and an \'etale sheaf of sets $\clF$ on $X$, according to the standard notation  we denote by $\clF_{\overline{x}}$ the stalk of $\clF$ at $x$ in the \'etale topology.  A gerbe  is an algebraic stack as in \cite[Definition 3.15]{LMBca}.

The main results in this paper are valid for banded $G$-gerbes over $X$ with $G$ a finite abelian group. For the sake of simplicity, in the main text we assume $G=\mu_r\subset \bbC^*$ is the cyclic group of $r$-th roots of unity. The case of general $G$ requires only notational changes. In genus $0$ this is spelled out in \cite[Appendix A]{AJT09-gerbes}. The higher genus case is similar and is left to the readers.

\subsection*{Acknowledgments} 
We thank D. Abramovich, A. Bayer,  K. Behrend,  B. Fantechi,  P. Johnson, A. Kresch, F. Nironi,  E. Sharpe,  Y. Ruan and A. Vistoli for valuable discussions. H.-H. T. is grateful to T. Coates, A. Corti, H. Iritani, and X. Tang for related collaborations. H.-H. T. is supported in part by NSF grant DMS-1047777.

\section{Orbifold Gromov-Witten theory}\label{OGWT-intro}
In this section we recall some basic set-up and well known facts about  orbifold Gromov-Witten theory that we will constantly use in what follows.

\subsection{Twisted stable maps}
We recall the definition of twisted curve here, see \cite{AGV01}, \cite{AGV06}, \cite{AV02} for more details.
\begin{defn}[\cite{AV02}, Definition 4.1.2]
A twisted nodal $n$-pointed curve over a scheme $S$ is a morphism $\mathcal{C}\to S$ together with $n$ closed substacks $\sigma_i\subset \mathcal{C}$ such that
\begin{itemize}
\item
$\mathcal{C}$ is a tame Deligne-Mumford stack, proper over $S$, and \'etale locally is a nodal curve over $S$;

\item
$\sigma_i\subset \mathcal{C}$ are disjoint closed substacks in the smooth locus of $\mathcal{C}\to S$;

\item
$\sigma_i\to S$ are \'etale gerbes;

\item
the map $\mathcal{C}\to C$ to the coarse moduli space $C$ is an isomorphism away from marked points and nodes.
\end{itemize}
\end{defn}

By definition the genus of a twisted curve $\mathcal{C}\to S$ is the genus of its coarse moduli space $C\to S$. 

Throughout this paper we will always assume that twisted curves are {\em balanced}, i.e. at any twisted node, the local group acts on the two branches by opposite characters.

 For more details on twisted curves that we will need in this paper we refer to \cite{AJT09-gerbes}, where the reader
will find an introduction to the equivalence between twisted curves and log twisted curves introduced in \cite{OLogCurv}.

Let $S$ be a noetherian  scheme and let  $\mathcal{X}/S$ be a proper  Deligne-Mumford stack over $S$ with projective coarse moduli space $X\to S$.
We fix an ample invertible sheaf $\clO_X(1)$ over $X$. 
 Let $\mathcal{K}_{g,n}(\mathcal{X}, \beta)$  be the fibered
category over $S$ which to any $S$-scheme $T$ associates the groupoid of the following data:
\begin{itemize}
\item
A twisted $n$-pointed curve $(\mathcal{C}/T, \{\sigma_{i}\})$ over $T$;
\item
A representable morphism $f : \mathcal{C}\to \mathcal{X}$ such that the induced morphism $\bar{f}: C\to  X$ between coarse moduli spaces is an $n$-pointed stable map of degree $\beta\in H^+_2(X, \mathbb{Z})$ (i.e. $\bar{f}_*[C]=\beta$).
\end{itemize}
According to \cite[Theorem 1.4.1]{AV02}, the fibered category $\mathcal{K}_{g,n}(\mathcal{X}, \beta)$ is a Deligne-Mumford stack  proper  over $S$. As discussed in \cite{AGV06}, there exists  evaluation maps: 
$$ev_{i}: \mathcal{K}_{g,n}(\mathcal{X},\beta)\to \bar{I}(\mathcal{X}), \quad 1\leq i\leq n$$ taking values in the {\em rigidified inertia stack} $\bar{I}(\mathcal{X})$ of $\mathcal{X}$. This map is obtained as follows. The rigidified inertia stack $\bar{I}(\clX)$ may be defined as the stack of cyclotomic gerbes in $\clX$, i.e. representable morphisms from cyclotomic gerbes to $\clX$. The evaluation map $ev_i$ is defined by associating to a twisted stable map $f: (\mathcal{C}/T, \{\sigma_{i}\})\to \mathcal{X}$ its restriction to the $i$-th marked gerbe,  
$$f|_{\sigma_i}: \sigma_i\to \mathcal{X},$$
which is an object of $\bar{I}(\mathcal{X})$. The rigidified inertia stack $\bar{I}(\clX)$ has an alternative description. Define the {\em inertia stack} of $\clX$ to be the fiber product over the diagonal: $$I\clX:= \clX\times_{\clX\times_S \clX}\clX.$$ 
By definition, objects of $I\clX$ are pairs $(x, g)$ where $x$ is an object of $\clX$ and $g$ is an element of the automorphism group of $x$. The rigidified inertia stack $\bar{I}(\clX)$ is obtained from $I\clX$ by applying the rigidification procedure (\cite{ACV}, \cite{AOV07}). More details can be found in e.g. \cite{AGV06}. 

\subsection{Virtual Fundamental Class}\label{vfc-subsection}
Let $\clX$ be a smooth proper DM stack with projective coarse moduli space. Let $X$ be a smooth projective variety. In general both $\clK_{g,n}(\clX,\beta)$ and $\barM_{g,n}(X,\beta)$ are neither smooth nor of pure dimension. Therefore they do not generally admit a fundamental cycle. In \cite{BFinc} the {\em virtual fundamental class} is introduced, which replaces the usual fundamental cycle. In order to construct it, an {\em  obstruction theory} is needed.\\
{\bf Condition ($\dagger$):}
We say that an object $L^\bullet$  of $D(\clO_{X_{\acute{e}t}})$ satisfies Condition $(\dagger)$
if
\begin{enumerate}
\item $h_i(L^\bullet )=0$ for all $i > 0$,
\item $h_i(L^\bullet )$ is coherent, for $i = 0, −1$.
\end{enumerate}
\begin{defn}\label{ob-thy-defn}
 An {\em obstruction theory} is a morphism in $D(\mathcal{X})$ $\phi: E^\bullet\to L_\mathcal{X}^\bullet$ where $ E^\bullet\in (\clO_{X_{\acute{e}t}})$ satisfies Condition ($\dagger$)  and  $L_\mathcal{X}$ is Illusie's cotangent complex (for extension to DM stacks see \cite{LMBca} and to Artin stacks \cite{Olss-Log-Cotg}). The morphism $\phi$ has to satisfy the following conditions:
\begin{enumerate}
 \item $\phi^0:h^0(E^\bullet)\to h^0(L_\mathcal{X}^\bullet)$ is an isomorphism;
\item $\phi^1:h^1(E^\bullet)\to h^1(L_\mathcal{X}^\bullet)$  is surjective;
\end{enumerate}
where $h^i$ denotes the $i$-th complex cohomology and $\phi^i$ is the morphism induced  by $\phi$ on cohomology sheaves.
\end{defn}
\begin{defn}
 An obstruction theory as in Defiinition \ref{ob-thy-defn} is called {\em perfect} if $E^{\bullet}$ is locally quasi-isomorphic to a complex of vector bundles. It is called {\em of perfect amplitude contained in $[-1,0]$} if it is locally quasi-isomorphic to a complex of vector bundles $F^{-1}\to F^0$.
\end{defn}

There is also a notion of  {\em relative obstruction theory} for  a morphism $\clX\to\mathfrak{Y}$ form a DM stack to a smooth Artin stack of pure dimension.  The definition is 
analogous to Definition \ref{ob-thy-defn} except that  the absolute tangent complex $L^\bullet_{\clX}$ is replaced by the relative cotangent complex $L_{\clX/\mathfrak{Y}}^\bullet $. 
Given any separated DM stack $\mathcal{X}$ endowed with a perfect  obstruction theory it is possible to associate to it a cycle in $A_*(\clX)_\bbQ$ of the expected dimension called   the {\em virtual fundamental class}. Such a cycle can be seen more conveniently as a bivariant class in $A^*(\clX\to S)\simeq A_*(\clX)$, $S=\Spec{\bbC}$.
Analogously, in the relative case we get  a class in $A_\bbQ^*(\clX\to\YY)\simeq A^*_\bbQ(\clX)$ (see \cite[Ch. 17]{Fu} for the theory of bivariant classes).
The construction works as follows. Let us observe first that   any DM type morphism of Artin stacks $\mathfrak{f}: \mathfrak{X}\to\mathfrak{Y}$ defines  a cone stack $\mathfrak{C}_{\mathfrak{X}/\mathfrak{Y}}$  called  {\em intrinsic relative normal cone}.  For any such a morphism one can find a diagram as follows \cite{Kre99}
\begin{eqnarray}
 \xymatrix{
U\ar[d]^{\acute{e}t} \ar@{^{(}->}[r] & M\ar[d]^{sm}\\
\mathfrak{X}\ar[r]^{\mathfrak{f}} & \mathfrak{Y}
}
\end{eqnarray}
where the upper row is a closed immersion, the arrow on the left is \'etale and the arrow on the right is smooth. This is  called in \cite{BFinc} a {\em local embedding}. Such a diagram defines a cone stack $[C_{U/M}/T_{M/\mathfrak{Y}}]$ which is shown to be in fact independent on the chosen local embedding of $\mathfrak{X}$. Moreover, the cone stacks constructed locally glue to yield a global cone stack.
\begin{defn}
The {\em intrinsic relative normal cone}  $\mathfrak{C}_{\mathfrak{X}/\mathfrak{Y}}$  is the unique cone stack  such that
$\mathfrak{C}_{\mathfrak{X}/\mathfrak{Y}}|_{U}\simeq [C_{U/M}/T_{M/\mathfrak{Y}}]$.
\end{defn}
Let $\clX\to \mathfrak{Y}$ be a morphism from a DM stack to a smooth Artin stack of pure dimension. Let $E^{\bullet}$ be a perfect obstruction theory of amplitude contained in $[-1,0]$ for $\clX$ relative to $\YY$. Let $E_\bullet:=(E_{fl}^\vee)^{\bullet}$, where the subscript $_{fl}$ denotes the derived pullback to the big fppf site. The quasi-isomorphism class of such an obstruction theory determines a  Picard stack $h^1/h^0((E_{fl}^\bullet)^\vee)\simeq [E_1/E_0]$ (the stack theoretic quotient). Note that if  $E^{\bullet}$ admits a {\em global resolution}, namely it is globally isomorphic to a complex of vector bundles $F^{-1}\to F^0$, then the associated Picard stack is simply $[F_1/F_0]$. The conditions in Definition \ref{ob-thy-defn} ensure that there is a closed embedding 
$\mathfrak{C}_{\mathfrak{X}/\mathfrak{Y}}\hookrightarrow [E_1/E_0]$.
 If $E^{\bullet}$ admits a  global resolution by  pulling-back $\mathfrak{C}_{\mathfrak{X}/\mathfrak{Y}}$ along $F_1\to [F_1/F_0]$ we get a closed substack  $C\hookrightarrow F_1$. The virtual fundamental class is obtained as the intersection of $C$ with the zero section of $F_1$.\footnote{In \cite{Kre99} intersection theory for Artin stacks was developed. Using such a theory it is possible to intersect $\mathfrak{C}_{\mathfrak{X}/\mathfrak{Y}}$ with the zero section of $[E_1/E_0]$ without  needing to assume  that the obstruction theory admits a global resolution.} 
 
The  stacks $\clK_{g,n}(\clX,\beta)$ and  $\barM_{g,n}(X,\beta)$ admit perfect  obstruction theories of amplitude contained in $[-1,0]$  relative to 
  $\MM_{g,n}^{tw}$ and $\MM_{g,n}$.  
In \cite{Beh97GW} it is shown that a perfect  obstruction theory relative to $\MM_{g,n}$  for     $\overline{M}_{g,n}(X,\beta)$ is given by 
\begin{equation*}
E^{\bullet}:=R\pi_*(f^*\Omega_X\otimes \omega_{\pi})\to  L_{ \overline{M}_{g,n}(X,\beta)/\mathfrak{M}_{g,n}},
\end{equation*}
where 
\begin{equation*}
\xymatrix{
C\ar[d]^{\pi}\ar[r]^{f} & X\\
\overline{M}_{g,n}(X,\beta)
}
\end{equation*}
is the universal stable map. Consider the universal twisted stable map to the gerbe $\clG$:
\begin{equation*}
\xymatrix{
\clC\ar[d]^{\tilde{\pi}}\ar[r]^{\tilde{f}}& \clG\\
\clK_{g,n}(\clG,\beta).
}
\end{equation*}
 According to \cite{AGV06} the moduli stack
$\clK_{g,n}(\clG,\beta)$ has a perfect obstruction theory $\tilde{E}^\bullet$ relative to $\MM^{tw}_{g,n}$ analogously constructed:
\begin{equation*}
\tilde{E}^{\bullet}:=R\tilde{\pi}_*(\tilde{f}^*\Omega_{\clG}\otimes \omega_{\tilde{\pi}})\to  L_{ \clK_{g,n}(\clG,\beta)/\mathfrak{M}^{tw}_{g,n}}.
\end{equation*}
The complex $\tilde{E}^{\bullet}$  turns out to be quasi-isomorphic to the pullback of $E^{\bullet}$ as an object in $\clD_{Coh}( \clK_{g,n}(\clG,\beta))$.

\subsection{Gromov-Witten invariants and Descendant Potential}
In this paper we  will work with cohomological GW-invariants. For us  $[\clK_{g,n}(\clX,\beta)]^{vir}$ and  $[\barM_{g,n}(X,\beta)]^{vir}$ will be classes in $H_*(\clK_{g,n}(\clX,\beta),\bbQ)$ and in $H_*(\barM_{g,n}(X,\beta),\bbQ)$. Let $\gamma_1,..,\gamma_n$ be homogeneous elements in $H^*(I\clX,\bbQ)$. Then   $n$-points genus $g$ orbifold GW-invariants are defined\footnote{Here we ignore the subtlety involving $I\clX$ and its rigidification.} as
\begin{eqnarray}\label{GWI-def-eqn}
 \langle \gamma_1,..,\gamma_n\rangle^{\clX}_{g,\beta}=\int_{[\clK_{g,n}(\clX,\beta)]^{vir}} \prod_{i=1}^n ev_i^*(\gamma_i).
\end{eqnarray}
The vector space $H^*(I\clX,\bbQ)$ can be made into a graded vector space in the following way:
\begin{eqnarray}
 H^*(I\clX,\bbQ)=\oplus_{\Omega} H^{a-age(\Omega)}(\Omega,\bbQ),
\end{eqnarray}
where the direct sum is over the connected components $\Omega\subseteq I\clX$ and $age(-)$ is a locally constant function on the inertia stack defined in \cite[Definition 7.1.1]{AGV06}. The genus $g$ {\em descendant potential} is defined as follows. Let $\clL_1,..,\clL_n$ be the tautological line bundles over $\clK_{g,n}(\clG,\beta)$, defined as the line bundles whose fiber at a point is the cotangent space of the corresponding coarse moduli space curve. For any $i=1,..,n$ let $\psi_i=c_1(\clL_i)$. Note that with this definition the classes $\psi_i$ are a pullback of the analogous classes in $H_*(\barM_{g,n}(X,\beta), \bbQ)$. Let $\alpha_1,..,\alpha_m$ be an additive basis of $H^*(\clX,\bbQ)$. Let   $t_{i,j}$, $i=1,..,m$, $j\geq 0$ be  supercommuting variables such that $\mbox{deg } t_{i,j}=\mbox{deg }\alpha_i$. We put $\gamma=\sum_{d=0}^\infty\sum_{i=1}^m t_{i,j}\tau_j\alpha_i$, where $\tau_j\alpha_i$ is defined as $\alpha_i\psi_*^j$ where $*=k$ if $\tau_j\alpha_i$ is inserted at the $k$-th place in (\ref{GWI-def-eqn}). Then the genus $g$ descendant potential is the series with coefficients in $\bbC[[t_{i,j}]]$ defined as follows
\begin{eqnarray}\label{orb-desc-potential-eq}
 \clF^g(\gamma,Q)=\sum_{n\geq 0}\sum_{\beta\in H_2(X,\bbZ)} \frac{Q^\beta}{n!}\langle\gamma^n\rangle=\nonumber\\
\sum_{n\geq 0}\sum_{\beta\in H_2^+(X,\bbZ)}\sum_{\stackrel{i_1,..,i_n}{j_1,..,j_n}} \frac{Q^{\beta}}{n!}\prod_{k=1}^nt_{i_k,j_k}\langle \prod_{k=1}^n \alpha_{i_k}\psi_k^{j_k}\rangle_{\clX,\beta}.
\end{eqnarray}
The genus $g$ descendant potential is an exponential generating function for genus $g$ descendant orbifold Gromov-Witten invariants.

\section{Twisted curves}\label{pre}
In this section we present a few  facts about prestable curves and twisted curves.
\subsection{The Picard group of prestable curves}
Let $C$ be a genus $g$ smooth  curve.  Let $\Pic{C}^0$ denote the subgroup of degree $0$ line bundles.
 The extension
\begin{eqnarray} \label{can_lift_deg_eq}
1\to \Pic{C}^0\to \Pic{C}\to \bbZ\to 1
\end{eqnarray}
is a semidirect product.  Moreover  $\Pic^0{C}$ is a divisible group and for any $r\in\bbN$ the kernel of the multiplication by $r$  is 
a subgroup  of rank $r^{2g}$ \cite{Mumf-Ab}.
Therefore, for any line bundle $L\in\Pic{C}$ such that $\mbox{deg}\ L$ is multiple of $r$,
there are $r^{2g}$ line bundles $N\in \Pic{C}$ such that $N^{\otimes r}\simeq L$. 
If $C$ is prestable, which means that it only admits ordinary double points as singularities,
then  $\Pic^0{C}$  is a semi-abelian variety (see e.g. \cite{LiuQ}) as described in the following Lemma.

\begin{lem}\label{Pic_div_group_sing_lem}
 Let $C$ be a genus $g$ nodal prestable curve. Then  there is an exact sequence
\begin{eqnarray}
 1\to T\to \Pic{C}\to \Pic{C^\nu}\to 1,
\end{eqnarray}
where $T$ is an algebraic torus and $C^\nu$ denotes the normalization of $C$. Therefore $\Pic{C}$ is a semi-abelian variety. 
\end{lem}
\begin{rmk}
The algebraic torus $T$ in the Lemma above is isomorphic to $\bbG_m^s$ where $s$ is the number of double points. It can be thought of as parametrizing descent data for line bundles. Indeed it is a well known fact  (see e.g. \cite{HaMo_ModCurv}) that the category of line bundles  over $C$ is equivalent to the category of line bundles over $C^\nu$ endowed with an identification of the fibers over the preimages of  the nodes. Isomorphisms 
of 1-dimensional $k$-vector spaces are given by $\bbG_{m,k}$. 
\end{rmk}
In the following it will be convenient to   make use of {\em dual graphs}  (introduced in \cite{BehMan}, Definition 1.1 and Definition 1.5).
Given a prestable curve the associated   dual graph encodes 
its  topological type. 
\begin{defn}
A {\em graph} $\tau$  is a quadruple $(F_\tau , V_\tau, j_\tau , \partial_\tau )$, where $F_\tau$ and $V_\tau$  are finite sets, $\partial_\tau: F_\tau  \to V_\tau$  is a map and $j_\tau : F_\tau \to F_\tau$ an involution. We call $F_\tau$ the set of flags,$V_\tau$    the set of vertices, $S_\tau = \{f ∈ F_\tau | j_\tau f = f \}$ the set of tails and $E_\tau = \{\{f_1 , f_2 \} \subset F_\tau | f_1 = j_\tau f_1 \}$ the set of edges of $\tau$. For $v \in V_\tau$  let $F_\tau (v) = \partial_\tau (v)$ and $\abs{v} = \# F_\tau (v)$, the valence of $v$.
\end{defn}
\begin{defn}\label{modular-graph-def}
A {\em modular graph} is a graph $\tau$ endowed with a map $g_\tau:V_\tau\to \bbZ_{\geq 0}$; $v\mapsto g(v)$. The number $g(v)$ is called the {\em genus} of $v$. 
\end{defn}
\begin{notation}
Let $\tau$ be a modular graph.
 We denote by $b_1(\tau)$ the first Betti number of  $\tau$ defined as
$b_1(\tau)=1-\abs{V_\tau} +\abs{E_\tau}$
\end{notation}

The following lemma characterizes the torsion subgroups of the  Picard group  of a  prestable curve. 
 \begin{lem}\label{r-tors-Picard-card-lem}
 Let $C$ be a nodal  prestable curve of dual graph $\tau$. Let $r\in\bbN$. Let $\Pic{C}[r]$ be the $r$-torsion part of the Picard group.  Then 
\begin{eqnarray}
 \abs{\Pic{C}[r]}= 2g(\tau)-b_1(\tau),
\end{eqnarray}
where  $b_1(\tau)=1 - \# V_\tau + \# E_\tau$ is the first Betti number of the graph $\tau$.
\end{lem}

\subsection{The Picard group of twisted curves}
As for ordinary prestable curves in order to describe the Picard group of a twisted curve    along with its   torsion subgroups    we  will make use of  dual graphs encoding the
 topological type   of the curve and  the isotropy groups of its special points.
In  \cite{AJT09-prod} we introduced {\em gerby modular graphs} generalizing Definition \ref{modular-graph-def}.
They will allows us  to  suitably label strata and irreducible components of   stacks parametrizing twisted curves and twisted stable maps and will  make much easier our notation in  Section \ref{proof-push-sect}.

\begin{notation}\label{labeling_set}
Let $\mho$ be a finite set and $\mathfrak{o}: \mho\to \bbN$, $\mathfrak{a}: \mho\to \bbQ$ two set maps.
\end{notation}
The following example is important. 
\begin{example}\label{labeling_set_from_stack}
Let $\clX$ be a smooth proper Deligne-Mumford stack. A triple $(\mho, \mathfrak{o}, \mathfrak{a})$ is obtained as follows. Let $\mho=\mho(\clX)$ be the set of connected components of the rigidified cyclotomic inertia stack of $\clX$. Let $\mathfrak{o}$ be the set map  $\mho(\clX) \to \bbN$ taking $U\in\mho(\clX)$ to the integer $r$ such that  $U\subseteq \overline{I}_{\mu_r}(\clX)$. Let  $\mathfrak{a}: \mho(\clX)\to\bbQ$ be the set map taking $U$ to its age. 
\end{example}
\begin{defn}[Gerby dual graph]\label{gerby-graph-def}
Let $(\mho, \mathfrak{o}, \mathfrak{a})$ be as in {Notation \ref{labeling_set}}. A {\em gerby modular graph} $\ttau$ associated to $(\mho, \mathfrak{o}, \mathfrak{a})$ is the data of an underlying  modular graph $\tau$ with a map $\mathfrak{g}:F_\ttau\to \mho$ such that $\mathfrak{g}(f)= \mathfrak{g}(f')$ whenever  the flags $f$, $f'$ form an edge  $\{f,f'\}\in E_\ttau$. We define $\gamma:=\mathfrak{o}\circ\mathfrak{g}$. 

Let $A$ be a semigroup with indecomposable zero.  A {\em gerby $A$-graph} is a gerby modular graph $\ttau$  whose underlying modular graph $\tau$ is endowed with an $A$-structure (i.e. a map $V_\tau\to A$).
\end{defn}

\begin{notation}\label{gerby-graph-notation}
Let  $\ttau$ be a gerby dual graph over  $\tau$. We establish the following conventional notations,
\begin{itemize}
 \item For any  edge $e_j\in E_\ttau$ we put $r_j:=\gamma(e_j)$. If $\clC$ is a twisted curve with gerby dual graph $\ttau$, 
we denote by  $\Gamma_{e_j}$  the automorphisms group of a  twisted node of $\clC$  
corresponding to $e_j$. Such a group is cyclic of order $r_j$;
\item for any flag $p_i\in F_\ttau/E_\ttau$ we put $b_i=\gamma(b_i)$.  If $\clC$ is a twisted curve with gerby dual graph $\ttau$, We denote by  $\Gamma_{p_i}$  the automorphisms group of a  stacky point   of $\clC$  
corresponding to $p_i$. Such a group is cyclic of order $b_i$;
\end{itemize}
\end{notation}

\begin{lem}
 Let $\clC$ be a twisted curve over $\Spec{\bbC}$ of gerby dual graph $\ttau$.  Let $\nu:\clC^\nu\to \clC$ be its normalization. 
Then there is an exact sequence
\begin{eqnarray}
 1\to \prod_{i=1}^s \bbC^*\to \Pic{\clC}\to \Pic{\clC^\nu}\to \prod_{i=1}^s \Pic{B\Gamma_{e_i}}\to 1
\end{eqnarray}
where $s=\abs{E_\ttau}$ and $B\Gamma_{e_i}$ are the gerbes corresponding to twisted nodes.
\end{lem}

\begin{pf}
 It follows from the normalization exact sequence
\begin{eqnarray}
 1\to \clO_\clC\to \nu_*\clO_{\clC^\nu}\to \oplus_{i=1}^s \clO_{B\Gamma_{e_i}}\to 1
\end{eqnarray}
by taking the invertible elements and by computing the long exact cohomology sequence.
\end{pf}

\begin{lem}\label{Pic-twisted-curve-lem}
 Let $\clC$ be a twisted curve over $\Spec{\bbC}$ with coarse moduli space $C$.  Let $\ttau$ be the dual graph of $\clC$.  With reference to 
Notation \ref{gerby-graph-notation},   there is an  exact sequence
  \begin{eqnarray}\label{sing-curve-pic-exact-seq}
 1 \to \Pic{C}\to \Pic{\clC}\to (\bigoplus_{e_j\in E_\ttau} \bbZ/r_j\bbZ\bigoplus_{p_i\in F_\ttau/E_\ttau} \bbZ/b_i\bbZ)\to 1.
\end{eqnarray}
\end{lem}

\begin{pf}
Consider the exact sequence of complexes 
\begin{eqnarray} 
 1\to \pi_*\bbG_m\to R\pi_*\bbG_m\to \frac{R\pi_*\bbG_m}{\pi_*\bbG_m}\to 1.
\end{eqnarray}
The Hypercohomology long exact sequence  yields
\begin{eqnarray}\label{Pic-twisted-curve-eq}
 1 \to H^1(C,\bbG_m)\stackrel{\pi^*}{\longrightarrow} H^1(\clC,\bbG_m)\stackrel{res}{\longrightarrow} H^0(C,R^1\pi_*\bbG_m)\to 1,
\end{eqnarray}
where the arrow on the right is surjective because of Tsen's theorem (cf. Section \ref{brauer_gp}).  
By a result proven in \cite[Proposition A.0.1]{ACV}  there is a canonical  isomorphism
\begin{eqnarray}\label{Pic-clC-mod-C-eq}
 H^0(C,R^1\pi_*\bbG_m)\simeq \prod_{e_j\in E_\ttau} H^1(\Gamma_{e_j}, \bbG_m)\times\prod_{p_i\in F_\ttau/E_\ttau} H^1(\Gamma_{p_i}, \bbG_m),
\end{eqnarray}
where the groups above are group cohomology for the trivial action of $\Gamma_{e_j}$, $\Gamma_{p_i}$ on $\bbG_m$. It is a standard result
that there are isomorphisms $H^1(\Gamma_{e_j}, \bbG_m)\simeq\prod_j\bbZ/r_j\bbZ$, $H^1(\Gamma_{p_i}, \bbG_m)\simeq \prod_i\bbZ/b_i\bbZ$.
 The map $res$ in (\ref{Pic-twisted-curve-eq})   is given by restricting   line bundles to the twisted points of $\clC$. 
\end{pf}

\begin{rmk}\label{multple-rmk-twisted-pic}
Let $\clC$ be a twisted curve with gerby dual graph $\ttau$. Let $\#E_\ttau=s$, $\#(F_\ttau\setminus E_\ttau)=n$. There is a commutative diagram
\begin{eqnarray}
 \xymatrix{
1\ar[r]& \Pic{C}\ar[r]\ar[d]& \Pic{\clC}\ar[d]^{\nu^*} \ar[r]& \oplus_{i=1}^n \bbZ/b_i\bigoplus \oplus_{k=1}^s \bbZ/r_k\ar[d]^{\mathbb{I}\bigoplus \Delta^{-1}}\ar[r] & 1 \\
1\ar[r] & \Pic{C^\nu} \ar[r] &  \Pic{\clC^\nu}\ar[r] \ar[d]& \oplus_{i=1}^{n}\bbZ/b_i\bigoplus \oplus_{k=1}^s \bbZ/r_k^{\oplus 2}\ar[ld] \ar[r] & 1\\
& &  \oplus_{k=1}^{s}\bbZ/r_k \ar[d] & &  \\
& & 1 &  & 
}
\end{eqnarray}
where the map $\mathbb{I}\bigoplus \Delta^{-1}$ is the identity of the first factor and the anti-diagonal on the second factor. The diagonal arrow is its cokernel. This corresponds to the fact
that a line bundle which is non trivial when restricted to a twisted node over $\clC$
pulls-back over $\clC^{\nu}$ to a line bundle involving opposite powers of the tautological line bundles
associated to the two preimages of the node. This is a consequence of the fact that nodes of twisted curves
are balanced.  We observe  that any two line bundles in $\Pic{\clC}$ pulling back to isomorphic line  bundles in $\Pic{\clC^\nu}$ differ by a degree zero line bundle because the normalization map is of degree 1. 
\end{rmk}

The following lemma gives a characterization of torsion subgroups for unmarked twisted curves.
\begin{lem}\label{r-tors-Picard-nonsep-lem}
 Let $\clC$ be an unmarked twisted curve with gerby  dual graph $\ttau$.
With reference to Notation \ref{gerby-graph-notation}  we  define   $l_e:=g.c.d.(r_e,r)$ $\forall e$.    Then there is a non split short exact sequence
\begin{eqnarray}
1\to \bbZ/r\bbZ^{2g(\ttau)-b_1(\ttau)}\to \Pic{\clC}[r]\to \bigoplus_{e\in E_{\ttau}}^{n.s.} \bbZ/l_e\bbZ\to 1,
\end{eqnarray}
where the direct sum $\bigoplus_{e\in E_{\ttau}}^{n.s.}$  is taken over the non separating nodes. 
\end{lem}
\begin{pf}
Note that $\bbZ/r\bbZ^{2g(\ttau)-b_1(\ttau)}=\bbZ/r\bbZ^{2g(\tau)-b_1(\tau)}\simeq\Pic{C}[r]$. We will see that line bundles which are non-trivial when restricted to non-separating nodes  contribute to $\Pic{\clC}[r]\simeq H^1(\clC,\mu_r)$, where we use the equivalence between the category of $r$-torsion line bundles and the category of $\mu_r$-principal bundles. Consider the exact sequence of complexes 
\begin{eqnarray}
1\to \pi_*\mu_r\to R\pi_*\mu_r\to \frac{R\pi_*\mu_r}{\pi_*\mu_r}\to 1.
\end{eqnarray}
By taking the Hypercohomology long exact sequence we get  the short exact sequence  
\begin{eqnarray}\label{ex-seq-mu-r}
 1\to H^1(C,\mu_r)\to H^1(\clC,\mu_r)\to \mbox{Ker}(H^0(C,R^1\pi_*\mu_r)\stackrel{\delta}{\to} H^2(C,\mu_r))\to 1.
\end{eqnarray}
The morphism $H^1(\clC,\mu_r)\to H^0(C,R^1\pi_*\mu_r)$ can be seen as the restriction map, which takes an $r$-torsion line bundle to its restriction to the stacky points.
The boundary morphism $\delta$ is explicitly described in \cite{Ch06}. Let $\vec{\alpha}=(\alpha_1,..,\alpha_m)$ be a class in  $H^0(C,R^1\pi_*\mu_r)\simeq \prod_{e\in E_\ttau} H^1(B\Gamma_e, \mu_r)\simeq \Pic({\coprod_{e\in E_\ttau} B\Gamma_e})$. Let $\clL_{\vec{\alpha}}$ be any  line bundle  over $\clC$ which restricts to a line bundle of class $\vec{\alpha}$ over $\coprod_{e} B\Gamma_e$. Then $\delta(\vec{\alpha})$ is the isomorphism class of the gerbe over $C$ induced by the $r$-th root of $\clL_{\vec{\alpha}}^{\otimes r}$. Note that there  is a commutative diagram (cfr. \cite[Lemma 3.2.19]{Ch06})
\begin{eqnarray}
  \xymatrix{
 H^1(C^\nu, R^1\pi^\nu_* \mu_r)  \ar[r]& H^2(C^\nu,\mu_r) \\
 H^1(C, R^1\pi_* \mu_r) \ar[u]_{\nu^*} \ar[r] & H^2(C,\mu_r)\ar[u]^{\wr},
 }
 \end{eqnarray}
where $C^\nu$  denotes the normalization of $C$. Using the notations in Remark \ref{multple-rmk-twisted-pic}, the above diagram becomes
\begin{eqnarray}
 \xymatrix{
\oplus_{k=1}^s {\bbZ/r_k}^{\oplus 2}\ar[r] & \oplus_{v\in \tau}(\bbZ/r)_v\ar[d]\\
\oplus_{k=1}^s \bbZ/r_k\ar[r]\ar[u]_{\Delta^{-1}} & \oplus_{v\in \tau}(\bbZ/r)_v\ar[u]^{\sim}
}
\end{eqnarray}
where for any $v$, $(\bbZ/r)_v=\bbZ/r$. The maps $\oplus_{k=1}^s \bbZ/rk^{\oplus 2}\to\oplus_{v\in \tau}(\bbZ/r)_v$ take an element $(m_1,...,m_{2s})$ in the image of $\Delta^{-1}$ to $\sum_{e_k\in E_v} (r/r_k) m_k$. Choose a set of line bundles whose images in $\Pic{\clC}/\Pic{C}$ generate it. Assume moreover that for any such line bundle there exists a  node $e$ in $\clC$ such that the restriction of this line bundle to the complement of $e$ in $\clC$ is trivial. The class of each line bundle of the above set is $(0,..,\alpha,..,0)$ with $\alpha$ in the $k$-th position. If $e_k$ is a non separating node, then $\alpha$ is sent to zero, because the preimages of the node in the normalization belong to the same irreducible component. If $e_k$ is separating, its preimages in $\clC^\nu$ belong to two different irreducible components, hence the image of $\alpha$ is $(\frak{r}{r_k}\alpha, -\frak{r}{r_k}\alpha)\in (\bbZ/r)_v\oplus(\bbZ/r)_{v'}$ where $v$, $v'$ are the vertices associated to the irreducible components of $\clC^\nu$ containing the preimages of the node. Therefore $\mbox{Ker}(\delta)$ is generated by the subset of line bundles whose restriction to non-separating nodes is non trivial. The claim follows from the fact that for any node $e$, $H^1(B\Gamma_{e},\mu_r)\simeq \bbZ/l_e$.  
\end{pf}

For any twisted curve $\clC$ we will fix a set of (isomorphisms classes of) line bundles in $\Pic{\clC}$ lifting the standard basis of $\Pic{\clC}/\Pic{C}$. We start from the right factor group in (\ref{Pic-clC-mod-C-eq}) in Lemma \ref{Pic-twisted-curve-lem}. Recall that any smooth twisted curve is obtained as stack over its coarse moduli space by the construction known as {\em taking roots of line bundles } as shown in \cite{AGV06}, \cite{Cadm03}. This construction provides for  any   marked twisted point $p_i$   a {\em tautological  line bundle} whose fiber over $p_i$ is the standard representation of $\Gamma_{p_i}\simeq \mu_{b_i}$. We will denote such line bundles by $\clT_i$ for any marked point $p_i$. If $p_i$ is twisted, $\clT_i$ has non-trivial image in $\Pic{\clC}/\Pic{C}$and it is our choice for the lift of its image. We  come to the second factor in (\ref{Pic-clC-mod-C-eq}). Unlike twisted marked points, for  twisted nodes  there is no universal construction providing a canonical choice for line bundles  whose restriction to a twisted node is non trivial. However,  due to Lemma \ref{r-tors-Picard-nonsep-lem} we can choose a set of line bundles, that we will denote by  $\clQ_1,..,\clQ_s$, $s=\abs{E_\ttau}$, satisfying the following properties

\begin{property}\label{canonical-tw-lb-nodes-property}
\hfill
\begin{enumerate}
\item for any $e_i\in E_\ttau$, $i=1,..,s$, $\clQ|_{e_i}$ is the standard representation of $\Gamma_{e_i}\simeq \mu_{r_i}$;
 \item for any $i=1,..,s$, $\clQ_i^{r_i}\simeq \clO_\clC$, $r_i=\gamma(e_i)$;
\item for any $i=1,..,s$ $\clQ_i|_{\clC\setminus e_i}\simeq \clO_{\clC\setminus e_i}$.
\end{enumerate}
\end{property}
We remark that any two line bundles $\clQ_i$, $\clQ'_i$ satisfying the above conditions for some $i$
differ by an $r_i$-torsion line bundle.\label{canonical-tw-lb-nodes-page}

\subsection{The Brauer group}\label{brauer_gp}
The Brauer group of a  smooth algebraic curve over $\Spec{\bbC}$ is trivial.
This is a consequence of Tsen's theorem as explained e.g. in \cite[III \S 2 Example 2.22 Case (d)]{Mil80}. In fact by the same argument the Brauer group of a prestable curve is also trivial. It is not hard to see that the same result also holds for the Brauer group of a twisted curve.
Consider sequence (\ref{ex-seq-mu-r}) in  Lemma \ref{r-tors-Picard-nonsep-lem}.
The result follows by taking the long exact cohomology sequence and
by the result in \cite[Proposition A.0.1]{ACV}. Indeed for $\Gamma$ any cyclic group $H^2(\Gamma,\bbG_m)=0$. See also
\cite{FPoma10}. 
\begin{lem}\label{Brauer_curve_over_artinian_ring}
 Let $\clC_A\to \Spec{A}$ be a twisted curve over an artinian ring. Then  $H^2(\clC_A, \clO_{\clC_A}^*)=0$.
\end{lem}
\begin{pf} 
Let $I\subset A$ be the ideal of the closed point $0$. Let $\clC_0$ denote the pullback of the curve to the closed point. By flatness $I\cdot \clO_{\clC_A}$ is the ideal sheaf
of $\clC^0$. The result follows by taking the long exact cohomology sequence of the short exact sequence
\begin{eqnarray}
 1\to 1+ I\cdot \clO_{\clC_A}\to \clO_{\clC_A}^*\to \clO_{\clC_0}^*\to 1,
\end{eqnarray}
and from $H^2(\clC_A, \clO_{\clC_A})=0$.
\end{pf}

\begin{lem}\label{etale-loc-root-gerbe}
Let $\pi:\clC_S\to S$ be a twisted curve over a scheme $S$. Let $\clG\to\clC$ be a $\mu_r$-banded gerbe of class $\alpha\in H^2(C,\mu_r)$. Then \'etale locally over $S$, $\clG$ is isomorphic to a {\em gerbe of roots of a line bundle}.
\end{lem}
\begin{pf}
As explained e.g. in \cite{Liebl04}, $\clG$ is a root gerbe if and only if it admits a line bundle whose fibers carry an action of $\mu_r$ given by multiplication by $\chi(\alpha)$ where $\chi$ denotes the natural inclusion of the sheaof of $r$-th roots of unity in $\clO_\clC$. By standard limit arguments, such kind of sheaf exists over a geometric fiber $\clG_{\overline{s}}$ of $\pi$ and extends to an \'etale neighborhood of $\overline{s}$. 
\end{pf}

\section{Rephrasing the moduli problem}
 Recall that for the purpose of computing GW invariants, we 
can pretend that there are evaluation morphisms taking values on the usual inertia stack (rather than on the rigidified inertia stack)
$$ev_i: \clK_{g,n}(\clG,\beta)\to I\clG\simeq \clG\times_X G,$$ where the last canonical isomorphism
is due to the fact tha $\clG$ is a $G$-banded gerbe.
 Given an $n$-tuple $\vec{g}\in G^{\times n}$ we define
\begin{eqnarray}
\clK_{g,n}(\clG,\beta)^{\vec{g}}:=\cap_{i=1}^n ev_i^{-1}(I\clG_{g_i}), 
\end{eqnarray}
where $I\clG_{g_i}:=\clG\times \{g_i\}$.  
The substack  $\clK_{g,n}(\clG,\beta)^{\vec{g}}$ is either empty or an open and  closed component of  $\clK_{g,n}(\clG,\beta)$. It turns out that $\clK_{g,n}(\clG,\beta)^{\vec{g}}\neq\emptyset$ if and only if $\vec{g}$ is $\beta$-admissible in the sense of Definition  \ref{adm-vect-defn}. By definition $\clK_{g,n}(\clG,\beta)=\coprod_{\vec{g}}\clK_{g,n}(\clG,\beta)^{\vec{g}}$.  
 
Let $\clG\to X$ be a $\mu_r$-gerbe. The moduli problem of twisted stable maps to $\clG$ relative to the moduli problem of stable maps to $X$  is equivalent to the moduli problem of twisted stable maps to   gerbes  induced by pullback over the prestable curves admitting stable maps to $X$. This fact follows  easily from the universal property of the fiber product in the strict 2-category of algebraic stacks.  The datum of a twisted stable map  $[\tilde{f}:\clC\to \clG]$ is equivalent to the outer part of the following diagram:
\begin{eqnarray}\label{fiber_prod_diagram}
 \xymatrix{
\clC\ar[rd]^h\ar@/^1pc/[rrd]^{\tilde{f}} \ar@/_1pc/[rdd]_\pi & & \\
& \clG_C\ar[r]\ar[d] \ar@{}[rd]|{\square} & \clG\ar[d]\\
& C\ar[r]_f & X.
}
\end{eqnarray}
By the universal property the morphism $\tilde{f}$ induces a unique representable morphism $h:\clC\to \clG_C$, where $\clG_C$ is the pull-back gerbe.

In order to benefit from the above reformulation of the moduli problem we need to characterize the gerbes obtained by pullback to the prestable curves admitting stable maps to $X$. Note that  for a family of prestable curve $\pi: C\to S$, the Brauer group of a  geometric fiber is isomorphic to the stalk of the constructible sheaf  $R^2\pi_*\mu_r$ at the corresponding geometric point (see e.g. \cite{Mil80}, VI 2 Corollaries 2.3 and 2.5). In  particular the Brauer group is  not constant over $S$. The class of the pull-back gerbe $\clG_C$ in $H^2(C,\mu_r)$ induces a  global section of  $R^2\pi_*\mu_r$  whose evaluation at any geometric point $\overline{s}$ of $S$ is the class of the gerbe induced by pull-back over $C_{\overline{s}}$. By applying the relative trace map to the section of $R^2\pi_*\mu_r$ defined by $[\clG_C]\in H^2(C,\mu_r)$, we get a function over $S$ with values in $\bbZ/r$. A priori such a function
is only locally constant. In fact we will prove it is  constant. The element of $\bbZ/r$ obtained by evaluating this function has the following geometric meaning. The gerbe $\clG_{C_{\overline{s}}}$ induced over any  geometric curve $C_{\overline{s}}$ is isomorphic to the gerbe of $r$-th root of some line bundle.  Any two line bundles over $C_{\overline{s}}$ of the same degree induce isomorphic root gerbes. Therefore  we can associate to the gerbe $\clG_{C_{\overline{s}}}$  the degree of any line bundle inducing (a gerbe isomorphic to) it. If $C_{\overline{s}}$ is reducible, we consider the total degree of such a  line bundle. Moreover, we can consider componentwise the degree mod $r$, since the $r$-th root gerbe of a line bundle whose degree is a multiple of $r$ is trivial. In Lemma \ref{restr-map-lemma} and Proposition  \ref{Trace-prop} we describe  the restriction map of the gerbe $\clG$ along stable maps to $X$. It turns out that once the curve class $\beta=f_*[C]$ is fixed, the  total degree 
of any line bundle inducing a gerbe isomorphic to the pullback of $\clG$ to the geometric fibers of the domain curve  is constant over
$\barM_{g,n}(X,\beta)$ and only depends on $\beta$ and on the class of $\clG$. 

\begin{lem}\label{restr-map-lemma}
Let $f:C\to X$ be an object of $\barM_{g,n}(X,\beta)(\bbC)$. Let $\clG\to X$ be a gerbe of class $\alpha\in H^2(X,\mu_r)$. Then $$f^*\clG\simeq \sqrt[r]{L/C}$$ for any  $L\in \Pic{C}$ with  $\mbox{deg}\ L=\alpha^{an}\cap \overline{\beta}$, where $\alpha^{an}$ is
the image of $\alpha$ in $H^2(X^{an},\bbZ/r)$ and $\overline{\beta}$ is the image of $\beta$ in $H_2(X^{an},\bbZ/r)$. Here (and henceforth) $\sqrt[r]{L/C}$ denotes the stack of $r$-th roots of the line bundle $L$.
\end{lem}

\begin{pf} 
Since we work over $\bbC$ we have  a canonical isomorphism $(\bbZ/r)_{_X}\isomto (\mu_r)_{_X}$. Therefore from now on we will identify $\mu_r$ and its tensor powers with $\bbZ/r$.  Moreover, because $\bbZ/r$ is a finite group, for any $i\geq 0$ there are isomorphisms between \'etale and complex cohomology (\cite[III \S 3.3]{Mil80})
\begin{eqnarray}\label{isom-et-cpx-eq}
H^i(X,\bbZ/r) &  \to  & H^i(X^{an},\bbZ/r).\\
\delta \quad \quad  &  \mapsto &  \quad \quad  \delta^{an}\nonumber
\end{eqnarray}
Since $f$ is a proper map of proper schemes over $\bbC$, a restriction map for the \'etale cohomology 
$f^*:H^i(X,\bbZ/r) \to H^i(C,\bbZ/r)$ is defined. A pushforward map $f_*:H^i(C,\bbZ/r)\to H^{i-2c}(X,\bbZ/r) $, where
$c:=\mbox{dim}\, X-\mbox{dim}\, C$, is also defined by duality.
A class   $\eta\in H^i(C,\bbZ/r)$ is taken 
to the unique class $f_*\eta$ such that the condition
$$
Tr_C(f^*\xi\cap \eta)=Tr_X(\xi\cap f_*\eta).
$$
is satisfied  for any  $\xi\in H^{2-i}(X,\bbZ/r)$.
The operation $\mbox{Tr}( f^*(-) \cap 1 )$, where $1$ is the generator of $H^0(C,\bbZ/r)$, defines an element $f_*1$  of the dual space of $H^2(X,\bbZ/r)$.
Such an element is the same for any map of class $\beta$. 
 Indeed, by  the isomorphism (\ref{isom-et-cpx-eq})  we get 
\begin{eqnarray}\label{tr-eq}
\mbox{Tr}_C ( f^*\alpha\cap 1) =\mbox{Tr}_{C^{an}} ( f_{an}^*\ \alpha^{an} \cap \overline{[C]}^{P.D.})=\mbox{Tr}_{X^{an}}(\alpha^{an}\cap \overline{\beta}^{P.D.}),
\end{eqnarray}
where $\overline{[C]}$ and  $\overline{\beta}$ are the images of the fundamental class $[C]$, resp.  of $\beta$ , in
$H_2(C^{an},\bbZ/r)$, resp. in $H_2(X^{an},\bbZ/r)$, and P.D. means the Poincar\'e dual.
\end{pf}

\begin{rmk}
Note that when there is $L\in \Pic{X}$ such that $\clG\simeq \sqrt[r]{L/X}$, by the functoriality of the $r$-th root construction, the restriction map takes  $\alpha \in H^2(X,\mu_r)$ to the class of the gerbe root of a line bundle of degree  $c_1(L)\cap \beta$  modulo $r$. Here  $\beta$ is seen as usual as an element of $\mbox{Hom}_{\bbZ}(\Pic{X},\bbZ)$. 
\end{rmk}

\begin{prop}\label{Trace-prop}
Consider the universal diagram
\begin{eqnarray}
 \xymatrix{
\mathfrak{C}\ar[d]_{\mathfrak{\pi}} \ar[r]^{\mathfrak{f}} & X \\
\overline{M}_{g,n}(X,\beta),
}
\end{eqnarray}
where $\mathfrak{C}$ is the universal curve and $\mathfrak{f}$ is the universal stable map. Let $e: H^2(\mathfrak{C}, (\mu_r)_{_{\mathfrak{C}}})\to H^0({\overline{M}_{g,n}(X,\beta)},R^2\pi_*\mu_r)$ be the edge map obtained from the Leray's spectral sequence for the sheaf $(\mu_r)_{_{\mathfrak{C}}}$ and the morphism $\pi$. Then the map
$$
\mbox{Tr}\  (e (\mathfrak{f}^*\sqcup) ):  H^2(X,\mu_r) \to (\bbZ/r)_{_{\overline{M}_{g,n}(X,\beta)}} 
$$
determines a constant global section of  $(\bbZ/r)_{_{\overline{M}_{g,n}(X,\beta)}}$.  The evaluation of such a section at a geometric point $\overline{p}$ gives the total degree of any line bundle  whose $r$-th root gerbe is  isomorphic to the  gerbe  $\clG_{C_{\overline{p}}}$ obtained by pulling back $\clG$  to   $C_{\overline{p}}$. Moreover  such  a degree depends linearly on $\beta$ modulo $r$. 
\end{prop}

\begin{pf}
Let $(\mu_r)_{_{\mathfrak{C}}}$ be the sheaf of $r$-th roots of unity over $\mathfrak{C}$. The higher direct image sheaves  $R^i\pi_*(\mu_r)_{_{\mathfrak{C}}}$ are constructible sheaves.
From the low degree exact sequence associated to the Leray's spectral sequence for the sheaf $(\mu_r)_{_{\mathfrak{C}}}$ and the morphism $\pi$, we get the following  
morphism of cohomology groups
\begin{eqnarray}\label{Leray-moduli-seq}
H^2(\overline{M}_{g,n}(X,\beta),\mu_r)\stackrel{H^2(\pi^*)}{\longrightarrow}H^2(C,\mu_r)\stackrel{e}{\longrightarrow} H^0(\overline{M}_{g,n}(X,\beta), R^2\pi_*\mu_r).
\end{eqnarray}
Note that  we denoted by $\mu_r$ both the sheaf $(\mu_r)_{_{\mathfrak{C}}}$ and  $\pi_*(\mu_r)_{_{\mathfrak{C}}}\simeq (\mu_r)_{_{\overline{M}_{g,n}(X,\beta)}}$, where the last isomorphism holds because families of prestable curves are geometrically connected. Let $\alpha=[\clG]\in H^2(X,\mu_r)$. The pullback $\mathfrak{f}^*\alpha\in H^2(\mathfrak{C},(\mu_r)_{_{\mathfrak{C}}})$ induces via the edge map $e$ a global section of the constructible sheaf $R^2\pi_*(\mu_r)_{_\mathfrak{C}}$. Let  $\overline{p}:\Spec{\bbC}\to \overline{M}_{g,n}(X,\beta)$  be a geometric point. Let  $C_{\overline{p}}$ be the pre-stable curve over $\overline{p}$ and $\iota_{\overline{p}}: C_{\overline{p}}\to \mathfrak{C}$ the inclusion. Let $\alpha_{\overline{p}}$ be the image of $\mathfrak{f}^*\alpha$ along the map $H^2(\iota_{\overline{p}}^*): H^2(\mathfrak{C}, (\mu_r)_{_{\mathfrak{C}}})\to  H^2(C_{\overline{p}}, (\mu_r)_{_{C_{\overline{p}}}})$. By the universal property of the quotient (obviously   $H^2(\iota_{\overline{p}}^*)\circ H^2(\pi^*)=0$)  it follows from sequence (\ref{Leray-moduli-seq}) that there is a factorization of $H^2(\iota_{\overline{p}}^*)$ through the edge map $e$ as follows
\begin{eqnarray}
 \xymatrix{
H^2(\mathfrak{C},(\mu_r)_{_{\mathfrak{C}}})\ar[d]^{H^2(\iota_{\overline{p}}^*)} \ar[rr]^{e}& & H^0(\overline{M}_{g,n}(X,\beta), R^2\pi_*(\mu_r)_{_{\mathfrak{C}}})\ar[d]^{H^0(\overline{p}^*)} \\
H^2(C_{\overline{p}},(\mu_r)_{_{C_{\overline{p}}}})\ar[r]^{\sim} &  H^0(\overline{p},(R^2\pi_*(\mu_r)_{_{C_{\overline{p}}}})) &    \ar[l]^{\sim}  H^0({\overline{p}},\overline{p}^*(R^2\pi_*(\mu_r)_{_{\mathfrak{C}}})).
}
\end{eqnarray}
where the isomorphism on the right of the last row holds  by proper base change theorem. Therefore  the evaluation of $e(\mathfrak{f}^*\alpha)$ at  a geometric point $\overline{p}$ gives $\alpha_{\overline{p}}$ in $\overline{p}^*R^2\pi_*(\mu_r)_{_\mathfrak{C}}\simeq H^2(C_{\overline{p}},(\mu_r)_{_{C_{\overline{p}}}})$. Since the universal curve is representable relative to its base, there is a trace map constructed as in the case of schemes:
\begin{eqnarray}\label{trace-map-eq}
 R^2\pi_*(\mu_r)_{_{\mathfrak{C}}}\stackrel{\sim}{\leftarrow} R^2\pi_!(\mu_r)_{_{\mathfrak{C}^{sm}}}\stackrel{\mbox{Tr}}{\longrightarrow} \bbZ/r,
\end{eqnarray}
where $\mathfrak{C}^{sm}\subseteq \mathfrak{C}$ is the smooth locus of the universal curve. The natural morphism $R^2\pi_{!}\mu_{r,_{\mathfrak{C}^{sm}}}\to R^2\pi_*\mu_{r,_{\mathfrak{C}}}$  in (\ref{trace-map-eq}) is an isomorphism because $\mathfrak{C}\setminus \mathfrak{C}^{sm}$ is  finite over the base. A priori $\mbox{Tr}(e(\mathfrak{f}^*\alpha))$ is only locally constant. However, by Lemma \ref{restr-map-lemma} we conlcude that it is in fact constant, because its evaluation depends only on $\alpha$ and $\beta$.
\end{pf}

\begin{notation}
In what follows for the sake of simplicity we will denote by $\alpha\cap \beta$ the element $\alpha^{an}\cap \overline{\beta}\in\bbZ/r$ computed  in Lemma \ref{restr-map-lemma}. This element characterizes the gerbe $\clG_C$ pulled back to a prestable curve  $C$ of class $\beta$. In fact it is the total degree of any line bundle over $C$ whose $r$-th root is isomorphic to $\clG_C$.
\end{notation}

\begin{defn}\label{adm-vect-defn}
Given an identification $G\simeq \mu_r$, a {\em $\beta$-admissible vector}  is an $n$-tuple $\vec{g}=(\exp(\frac{2\pi \sqrt{-1} m_1}{b_1}),.., \exp(\frac{2\pi \sqrt{-1}m_n }{b_n}) )\in G^{\times n}$, with  $b_i|r$ and $(m_i,b_i)=1$ for $i=1,..,n$, and such that 
$$
\prod_{i=1}^n \exp(\frac{2\pi \sqrt{-1} m_i}{b_i}) = \exp(\frac{2\pi \sqrt{-1} k}{r}),
$$
where $k=\alpha\cap\beta$ mod $r$.
\end{defn}

\begin{lem}\label{repr-morph-lemm}
 Let $C$ be a prestable curve of dual graph $\tau$. Let $\clG\simeq \sqrt[r]{L/C}$ for some
$L\in\Pic{C}$. Let $\clC$ be a twisted curve of gerby graph $\ttau$ over $C$ admitting a representable morphism to $\clG$.
By the universal property of root gerbes, there exists a pair  $(\clN,\phi)$ such that 
$\phi:\clN^{\otimes r}\simeq L$. Let $n=\abs{F_\ttau\setminus E_\ttau}$. Then
\begin{eqnarray}
 \clN\simeq \otimes_{i=1}^n  \clT_i^{m_i}\otimes M,
\end{eqnarray}
where $(m_i,b_i)=1$, the restriction of $M$ to the marked points is trivial  and for any $e\in E_\ttau$ $M_{|_e}$ is a faithful representation of $\Gamma_e\simeq \mu_{r_e}$. 
\end{lem}
\begin{pf}
 Let us consider the normalization $\nu:\clC^\nu\to \clC$. This is a representable morphism.
Let $l=1,..,m$ be an index labeling the irreducible components of $\clC$. 
 For any $l$ there are induced morphisms $\clC_l^\nu\to \clC_l$, with $\clC_l^\nu$ a smooth twisted curve.
Let $n_l:= \abs{F_{\ttau(v_l)}\setminus E_{\ttau(v_l)} }$, $k_l:=\abs{E_{\ttau(v_l)}^{n.l.}}$, $s_l:=\abs{E_{\ttau(v_l)}^{loop}}$, where $\ttau(v_l)$ is the gerby graph of the irreducible component $\clC_l$.
Then
\begin{eqnarray}
 \clN|_{\clC_l^\nu}\simeq \otimes_{i=1}^{n_l}\clT_i^{m_i}\bigotimes\otimes_{j=1}^{k_l} \clT_j^{m_j}\bigotimes
\otimes_{t=1}^{s_l}\clT_t^{m_t} \otimes M'
\end{eqnarray}
where $M'$ is pelled-back from the coarse moduli space $C_l^\nu$.
Note that there is a bijection between marked points and separating nodes of $\clC_l$ and of its normalization.
On the contrary, for any non-separating node of $\clC_l$ there are two marked points in its preimage in $\clC_l^\nu$.
To study representability of $f^\nu_l:\clC_l^{\nu}\to \clG$ by  \cite[Lemma 4.4.3]{AV02}, it suffices to study the homomorphism 
\begin{equation}\label{hom_between_aut}
Aut(\sigma_i)\to Aut(f_l^\nu(\sigma_i)),
\end{equation}
 induced by $f^\nu_l$ on stack points. Here $\sigma_i$ can be either a stack point mapping to a smooth point
or to a node in $\clC_l$.
 By $\sigma_i$ we mean precisely a morphism 
$\tilde{h}_i: \Spec{K}\to \clC$ from  an algebraically closed field  $K$ to $\clC_l^\nu$ with image in   the special locus. 
 By the root construction description of $\clC^\nu_l$ (see e.g. \cite[Example 2.7]{Cadm03} and \cite[Section 4.2]{AGV06}), the stack point $\sigma_i$ is equivalent to the data $(h_i, \clN_i, \tau_i,\phi_i)$, where $h_i:\Spec{K}\to C$ with image $p_i$, $\clN_i$ is a line bundle over $\Spec{K}$,  $\phi_i:\clN_i^{\otimes r_i}\stackrel{\sim}{\to} h_i^*\clO(p_i)$, $\tau_i$ is a section of $\clN_i$ such that $\phi_i(\tau_i^{r_i})= h_i^* s_i$, hence $\tau_i=0$. The image $f_l^\nu(\sigma_i)$ is given by $\tilde{h}_i^*\nu^*\clN$ and $\tilde{h}_i^*\psi: \tilde{h}_i^*\nu^*\clN^{\otimes r}\simeq \tilde{h}_i^*\pi^*f^*L$. 
Note that $\tilde{h}_i^*\clT_i$ is naturally isomorphic to $\clN_i$. An automorphism  $\epsilon\in Aut(\sigma_i)\simeq \mu_{r_i}$ is mapped to $\epsilon^{m_i}\in Aut(f^\nu_l(\sigma_i)))\simeq \mu_r$ since   $\clN=\pi^*M'\otimes \bigotimes_{i=1}^n \clT_i^{m_i}$.
We conclude by observing that for any $\sigma_i$ in $\clC_l^\nu$, the restriction of the normalization map
is an isomorphism with its image. Hence in particular for any $e$ in $E_\ttau$ $M_{|_e}$ is isomorphic to $\clT_i^{m_i}$ over $\sigma_i$ for some $\sigma_i$ in $\clC^\nu$.
Therefore  $M_{|_e}$ is a faithful representation of $\Gamma_e\simeq Aut(\sigma_i)$ because $\clT_i^{m_i}$ is.
\end{pf}
\begin{lem}\label{tw-st-map-adm-vect-lem}
Let $[\tilde{f}: (\clC, \{\sigma_i\})\to \clG]\in \clK_{g,n}(\clG, \beta)^{\vec{g}}(\bbC)$ with $\clG\simeq\sqrt[r]{L/X}$
for some $L$ in $\Pic{X}$. Then $\tilde{f}$ is equivalent to a line bundle
\begin{eqnarray}
 \clN\simeq \otimes_{i=1}^n\clT_i^{m_i}\bigotimes M
\end{eqnarray}
where $M$ is an in Lemma \ref{repr-morph-lemm} and for $i=1,..,n$ $m_i$ are determined by $\vec{g}$.
\end{lem}
\begin{pf}
  By definition the morphism $\tilde{f}|_{\sigma_i}:B\mu_{r_i}\simeq \sigma_i\to \clG$ is equivalent to an injective  homomorphism  $$\mu_{r_i}\hookrightarrow \mu_r, \quad \exp(2\pi \sqrt{-1}/r_i)\mapsto g_i.$$
The argument in the proof of Lemma \ref{repr-morph-lemm}, applied to the irreducible component of $\clC$  containing $\sigma_i$, shows that we may write 
\begin{equation}\label{elements_in_adm_vector}
g_i=\exp(2\pi\sqrt{-1}\frac{m_i}{b_i}), \quad \text{with } 0\leq m_i< b_i, \text{ and } (m_i, b_i)=1.
\end{equation}
Furthermore, if $\clL^{1/r}$ is the universal $r$-th root of $\clL$ over $\clG$, then $\tilde{f}|_{\sigma_i}^*\clL^{1/r}$ is the $\mu_{b_i}$-representation on which the standard generator $\exp(2\pi\sqrt{-1}/b_i) \in \mu_{b_i}$ acts by multiplication by $\exp(2\pi\sqrt{-1}m_i/b_i)$. In other words 
\begin{equation}\label{age_of_universal_root_line_bundle}
\mbox{age}_{\sigma_i}(\tilde{f}^*\clL^{1/r})=\frac{m_i}{b_i}. 
\end{equation}
\end{pf}

\begin{lem}\label{separating_nodes_order_lem}
Let $[\tilde{f}:\clC\to \clG]$ be an object of $\clK_{g,n}(\clG,\beta)^{\vec{g}}(\bbC)$ with underlying map between coarse spaces being $[f: C\to X]$. Then the orders of the isotropy groups of marked points and separating nodes  are determined by $\vec{g}$.
\end{lem}
\begin{pf}
The proof proceeds by induction on the number of irreducible components of $\clC$.
If $\clC$ is irreducible the claim follows from Lemma \ref{tw-st-map-adm-vect-lem}. 
Let us prove the induction step. Let $\clC_1\subset \clC$  be an irreducible component containing only a separating node $e$.
Let $\clC_2=\overline{\clC\setminus \clC_1}$. Set $\beta_1=f_*[C_1]$. Then $\tilde{f}_{|_{\clC_1}}$ is equivalent
to a line bundle $\clN_1$  of degree $k_1/r$, where $k_1=c_1(L)\cap\beta_1$.
By Riemann-Roch for twisted curves 
\begin{eqnarray}\label{ord-node-def-eqn}
 \frac{k_1}{r}-\sum_{i=1}^n age_{p_i}(\clN_1)=age_e(\clN_1)+N =\frac{k_1}{r}-\sum_{i=1}^{n_1}\frac{m_i}{b_i}=\frac{m_e}{b_e}+N
\end{eqnarray}
where $N\in\bbN$.
Note that the  preimage of $e$ along the closed immersion $i:\clC_1\hookrightarrow \clC$ is a smooth point, that we denote
again by $e$. By taking the fractional part of (\ref{ord-node-def-eqn}) we get
\begin{eqnarray}
\langle\frac{k_1}{r}-\sum_{i=1}^{n_1}\frac{m_i}{b_i}\rangle=\frac{m_e}{b_e}\quad (m_e,b_e)=1.
\end{eqnarray}
We conclude by observing that the order of  the node $e$ is given by $b_e$.
\end{pf}
\begin{rmk}
 We observe that unlike separating nodes, $\vec{g}$ does not determine the order of non-separating nodes.
In Lemma \ref{separating_nodes_order_lem} we used Riemann-Roch componentwise in order to determine
the order of the inertia group of separating nodes. That works because the preimage of a separating node
along the inclusion of an irreducible compoents is a smooth point. Of course this is false for a  non-separating nodes.
Riemann-Roch formula for twisted curve does not take into account the ages of line bundles at nodes.
This is a consequence of the fact that twisted nodes are assumed to be balanced.
\end{rmk}

\begin{defn}\label{vecg-compatible-gerby-graph-def}
Let $\tau$ be a dual graph. Let $\vec{g}$ be an admissible vector. We say that  a dual graph $\ttau$ over $\tau$ is {\em $\vec{g}$-compatible} if the following conditions are satisfied:
\begin{enumerate}
\item[i)]  for all $p_i\in F_\ttau/E_\ttau$, we have $\gamma(p_i)=b_i$ where $b_i$ is  determined by $\vec{g}$ as in Definition  \ref{adm-vect-defn}.

\item[ii)] for all non-looping  $e_j\in E_\ttau$, we have $\gamma(e_j)=s_j$  where $s_j$ is  determined by $\vec{g}$ as in the proof of Lemma \ref{separating_nodes_order_lem}.
\end{enumerate}
\end{defn}

The definition of $\beta$-admissible vector is  useful to characterize twisted stable maps to arbitrary  banded gerbes. Let $\clG$ be a $\mu_r$-banded gerbe. Let $\alpha\in H^2(X,\mu_r)$ be the class of $\clG$. Let  $[f:C\to X]$ be    a geometric  point   in $\overline{M}_{g,n}(X,\beta)$.  Let $k$ be a lift to $\bbZ$ of $\alpha\cap\beta$. Let us  choose a degree $k$ line bundle  $L$ over $C$ and an isomorphism $\theta$ such that $\theta: f^*\clG\isomto \sqrt[r]{L/C}$. Let us consider the diagram
\begin{eqnarray}\label{adm-vect-diagram}
 \xymatrix @R=0.3pc{
\clC\ar@/^2pc/[rrr]^{\tilde{f}}\ar[rr]^{h}\ar[rdd]_{(\clN,\phi)} &  & \clG_C \ar[r]^{g}\ar[ddd]\ar@{}[rd]|{\square} & \clG\ar[ddd]\\
& \alpha\ \Rightarrow & & \\
& \sqrt[r]{L/C}\ar[ruu]_{\theta^{-1}} & \\
& & C \ar[r]_{f} & X.
}
\end{eqnarray}
We claim that the  twisted stable map $\tilde{f}:\clC\to \clG$ determines a $\beta$-admissible vector $\vec{g}$ which can be computed as  the $\beta$-admissible vector of $(\clN,\phi)$. Indeed it turns out that $\vec{g}$ is well defined in the sense that it does not depend on the choice of $L$ and $\theta$. First we know that $(\clN,\phi)$  determines some $\vec{g}$  as in Lemma \ref{tw-st-map-to-banded-gerbe}. Let  $\{p_i\}$ be the set of marked points in  $\clC$. Let us consider the restrictions $h|_{p_i}$ and $\tilde{f}|_{p_i}$. They define (up to the choice of a section) objects of the cyclotomic inertia stacks $I_\mu\clG_C$ and $I_\mu \clG$. We recall that  over $\bbC$ the cyclotomic inertia stack is  canonically isomorphic to the usual inertia stack.  By definition of banded gerbe, there are canonical isomorphisms  $I\clG\simeq \clG\times_X (\mu_r)_{_X}$ and  $I\clG_C\simeq \clG_C\times_C (\mu_r)_{_C}$, where $(\mu_r)_{_C} \simeq f^{-1}(\mu_r)_{_X}$. Indeed $\clG_C$ is canonically banded by the restriction of $(\mu_r)_{_X}$ to $C$ (\cite[IV 2 Corollary 2.2.4]{Gir}). There is a 2-commutative diagram
\begin{eqnarray}\label{inertia-banded-gerbes-diag}
 \xymatrix{
I\sqrt[r]{L/C}\ar[r]^-{\sim}\ar[d]\ar@{}[rd]|{\Rightarrow} &   \sqrt[r]{L/C}\times_C (\mu_r)_{_C}\ar[d]^{(\theta^{-1},id)}\\
I\clG_C\ar[d]  \ar@{}[rd]|{\Rightarrow} \ar[r]^-{\sim} & \clG_C\times_C (\mu_r)_{_C} \ar[d]^{(g,f_*)}\\
I\clG \ar[r]^-{\sim} & \clG\times_X (\mu_r)_{_X},
  }
\end{eqnarray}
where the arrows on the left are  the morphisms between the inertia stacks induced by  $\theta^{-1}$ and $g$  while $f_*$ denotes the morphism of schemes corresponding to the morphism of sheaves $(\mu_r)_{_X} \to  f_*(\mu_r)_{_C} $.  The restriction  $h|_{p_i}$ is an object of $I\sqrt[r]{L/C}$, therefore by (\ref{inertia-banded-gerbes-diag}) and  by the projection $\clG_C\times_C (\mu_r)_{_C}\to (\mu_r)_{_C}$ it defines an object of $\mu_r$. By the 1-commutativity of the upper part of  diagram (\ref{adm-vect-diagram}), we see that $\tilde{f}|_{p_i}$ is associated to the same element of $\mu_r$. Once we choose $\theta$ and $L$ we get a bijection between the possible representable morphisms $\clC\to\clG_C$ and the representable morphisms $\clC\to \sqrt[r]{L/C}$. Since $\theta$ is an isomorphism of banded gerbes, which is the same as to say that the upper square of  diagram (\ref{inertia-banded-gerbes-diag}) commutes, any two morphisms  $(\clN,\phi)$, $h$ corresponding to each other via $\theta$ define the same $n$-tuple of elements of $\mu_r$. Therefore we can associate to a  twisted stable map $[\tilde{f}: \clC\to \clG]$ over $[f:C\to X]$  the  $n$-tuple of elements of $\mu_r$  defined by the morphism $\clC \to \sqrt[r]{L/C}$  associated to $\tilde{f}$ by any choice of $\theta$ and $L$. Such an $n$-tuple by definition is a  $\beta$-admissible vector. 
We will sometimes call a twisted stable map with $\beta$-admissible vector $\vec{g}$ a  twisted stable map of {\em contact type} $\vec{g}$.

%

The above discussion allows to state the following Lemma.
\begin{lem}\label{tw-st-map-to-banded-gerbe}
 Let $\clG\simeq \clG$ be a $\mu_r$-banded gerbe. Let $[f:C\to X]$ be an object of  $\overline{M}_{g,n}(X,\beta)(\bbC)$. Let $\clC\to C$ be a twisted curve over $C$ such  that the order of all of its special points divides $r$. Then a representable morphism
$\clC\to f^*\clG$  lifting $f$ is equivalent to a pair $(\clN,\phi)$, where $\clN\in\Pic{\clC}$ and $\phi: \clN^r\to f^*\clL$ satisfy the following conditions
\begin{enumerate}
\item $\phi:\clN^r\isomto f^*L$;
 
\item $\clN\simeq  \otimes_{i=1}^n \clT_i^{m_i}\otimes \clM$, where $\clT_i$ are the tautological line bundles associated to the marked points and $\clM\in \Pic{\clC}$  such that, for any marked point $p_i$,  $\clM|_{p_i}\simeq \clO_{p_i}$;

\item $\vec{g}\in \mu_r^{\times n}$, defined as the $n$-tuple $(g_1,..,g_n)$ such that  $g_i=\exp(\frac{2\pi \sqrt{-1} m_i}{b_i})$, is $\beta$-admissible;

\item  for any twisted node  $e\in \clC$, $\clM|_{e}$ is a vector bundle carrying a faithful representation of the isotropy group of $e$. 
\end{enumerate}
Moreover if  $[\tilde{f}:\clC\to\clG]$ is in $\clK_{g,n}(\clG,\beta)^{\vec{g}}$, then 
the order of the inertia groups of marked points and  of separating nodes is determined. In oter words, the dual graph
$\ttau$ of $\clC$ is $\vec{g}$-compatible over the dual graph $\tau$ of its coarse moduli space $C$.
\end{lem}
\begin{pf}
 It follows from Lemmas \ref{repr-morph-lemm}, \ref{separating_nodes_order_lem} and from the discussion starting on page \pageref{adm-vect-diagram}
\end{pf}

\begin{prop}\label{geom-lift-prop}
 Let $\clG\to X$ be  a $\mu_r$-banded  gerbe. Let $[f:(C, p_1,...,p_n) \to X]$ be an object of $\overline{M}_{g,n}(X,\beta)(\bbC)$. Let $\tau$ be the dual graph of $C$. 
  Let   $\vec{g}=(g_1,...,g_n)\in \mu_r^{\times n}$  a $\beta$-admissible vector. 
Let $\ttau$ be a  $\vec{g}$-compatible gerby graph such that for any $e\in E_\ttau^{loop}$ $\gamma(e)|r$.
Let $\clC$ be the twisted curve over $C$ with dual gerby graph $\ttau$.
Then $\clC$ admits twisted stable maps to $\clG$ with associated $\beta$-admissible vector $\vec{g}$ lifting  $[f:C\to X]$.
Moreover, the number of non isomorphic twisted stable maps is
\begin{eqnarray}
 N(\ttau)=r^{2g-b_1(\tau)}\times \prod_{l=1}^{\abs{E_\ttau}} \phi(\gamma(e_l))
\end{eqnarray}
where $\phi(\ )$ denotes the Euler totient function.
\end{prop} 
\begin{pf}
We sketch the proof, details are left to the reader. We fix a line bundle $L$ in $\Pic{C}$ such  that $f^*\clG\simeq \sqrt[r]{L/C}$. Let $\clC$ be a twisted curve as above.
Then a twisted stable map $\clC\to\clG$ lifting $f$ corresponds to a pair $(\clN,\phi)$ such that $\phi:\clN^{\otimes r}\simeq L$. By Lemma \ref{tw-st-map-to-banded-gerbe}
we know that the admissible vector $\vec{g}$ fixes  the restriction of $\clN$ to $\oplus_{p_i\in F_\ttau/E_\ttau} B\mu_{b_i}\bigoplus \oplus_{e_j\in E_\ttau^{n.l.}} B\mu_{r_j}$. By definition of $\beta$-admissible vector and by discussion on page \pageref{canonical-tw-lb-nodes-page} it is possible to lift this image to a line bundle $\clN'$ which is an  $r$-th root of $L$ and it is trivial when restricted to the complement of marked points and separating nodes. In order to get a line bundle $\clN$ corresponding to a representable morphism (cfr Lemma \ref{repr-morph-lemm}) we need to tensor $\clN'$ with an  $r$-torsion line bundle whose restriction to non-separating nodes yields  a faithful representation of their inertia groups and which is trivial elsewhere. Such a line bundle can be obtained as a tensor product of  line bundles $\clQ_l$ satisfying Property \ref{canonical-tw-lb-nodes-property}. The condition on the faithfulness of the representations of the inertia groups of the nodes is equivalent to the requirement for each $e_l\in E_\ttau^{loop}$ the associated line bundle $\clQ_l$ appears raised to a tensor power coprime with $\gamma(e_l)$ (see again  Lemma \ref{repr-morph-lemm}). Moreover, $\clN$ can be further tensored by $r$-torsion line bundles pulled-back from $C$ to obtain non-isomorphic $r$-th roots of $L$. The number of non-isomorphic $(\clN,\phi)$ yielding non-isomorphic twisted stable maps to $\clG$ is computed as 
\begin{eqnarray}\label{number-of-lifts-ttau}
 N(\ttau)=r^{2g-b_1(\tau)}\times \prod_{l=1}^{\abs{E_\ttau}} \phi(\gamma(e_l)),
\end{eqnarray}
where $\phi(-)$ denotes the Euler totient function, whose evaluation on $n\in\bbN$ gives the number of integers less than $n$ and coprime with $n$.
\end{pf}

\section{The structure morphism $p$}
As sketched in the introduction, in order to compare orbifold Gromov-Witten invariants of $\clG$ with Gromov-Witten invariants of $X$ we need to prove a push-forward formula for the virtual fundamental class of $\clK_{g,n}(\clG,\beta)$. We consider the following diagram:
\begin{eqnarray}\label{VFC_push_diag}
\xymatrix{
\clK_{g,n}(\clG,\beta)^{\vec{g}}\ar[dr]\ar[r]^{t}\ar@/^2pc/[rr]^{p} & P_{g,n}\ar[r]^-{q'}\ar[d]\ar@{}[rd]|{\square} & \barM_{g,n}(X,\beta)\ar[d] \\
& \MM_{g,n}^{tw} \ar[r]_{q} & \MM_{g,n}.
}
\end{eqnarray}
 Here the morphism $q$ maps a  twisted curve  to the underlying prestable curve. The right vertical arrow is the forgetful morphism taking a stable map $[f:C\to X]$ to the prestable curve $C$.  The stack  $P_{g,n}$ is   defined as the fiber product.  There is a natural morphism $p: \clK_{g,n}(\clG,\beta)^{\vec{g}}\to \barM_{g,n}(X,\beta)$ associating to a twisted stable map $[\tilde{f}:\clC\to\clG]$ the underlying stable map  $[f:C\to X]$ between the coarse moduli spaces. Just like the case of $\barM_{g,n}(X,\beta)$ there is a natural forgetful morphism  $\clK_{g,n}(\clG,\beta)^{\vec{g}}\to  \MM_{g,n}^{tw}$ taking a twisted stable map $[\tilde{f}:\clC\to\clG]$ to the  twisted curve $\clC$. The square defined by the forgetful morphisms and the morphisms $q$ and $q'$ is commutative. Therefore by the universal property of the fiber product there is an induced  morphism  $t$.
 
We will prove that the pushforward map
$$p_*: H_*(\clK_{g,n}(\clG,\beta)^{\vec{g}},\bbQ) \to H_*(\overline{M}_{g,n}(X,\beta),\bbQ)$$
takes $[\clK_{g,n}(\clG,\beta)^{\vec{g}}]^{vir}$ to a multiple of $[\overline{M}_{g,n}(X,\beta)]^{vir}$, see Theorem \ref{push-thm} below. In order to prove this push-forward formula we will need the following Proposition, which says that the morphism $t$ in (\ref{VFC_push_diag}) is \'etale.  

\begin{prop}\label{t-etaleness-prop}
The  morphism $t$  in diagram (\ref{VFC_push_diag}) is \'etale.
\end{prop}
\begin{pf}
Since all our stacks are locally noetherian locally of finite type, to prove that  a morphism is \'etale we can use the infinitesimal lifting criterion and check it over square zero extensions of Artinian local rings (\cite[\S 17]{EGAIV_IV}). Let
\begin{eqnarray}\label{sq-zero-ext}
 1\to I \to B\to A\to 1
\end{eqnarray}
be a square zero extension of local Artinian rings. Assume that we have  the following 2-commutative diagram
\begin{eqnarray}
 \xymatrix{
\Spec{A}\ar@{_(->}[d]\ar[r]\ar@{}[rd]|{\Rightarrow} & \clK_{g,n}(\clG,\beta)\ar[d]\\
\Spec{B}\ar[r] & P_{g,n}.
}
\end{eqnarray}
This is the datum of  an object $(\clC_B, f_B: C_B\to X)$ of $P_{g,n}(B)$, an object $(\clC_A, \tilde{f}_A: \clC_A\to \clG)$ of $\clK_{g,n}(\clG,\beta)(A)$ and a pair of isomorphisms $\varphi\in \mbox{Mor}\  \barM_{g,n}(X,\beta)(A)$ and $\psi\in \mbox{Mor}\ \MM_{g,n}(A)$ making the following diagram commutative
 \begin{eqnarray}
  \xymatrix{
 \clC_A\ar[r] \ar[d] \ar[r]_{\sim}^{\psi}   & \clC_B|_A\ar[d]  \\
 C_A\ar[r]_{\sim}^{\varphi}\ar[d]_{f_A}  & C_B|_A  \ar[d]^{f_B|_A} \\
  X \ar[r]_{id} & X. 
 }
 \end{eqnarray}
Here $f_A$ denotes the morphism induced by $\tilde{f}_A$ by  passing to the coarse moduli spaces. The morphism $t$ is \'etale if and only if there exists a unique (up to isomorphism) arrow $\Spec{B}\to \clK_{g,n}(\clG,\beta)$ making the following diagram 2-commute
\begin{eqnarray}
 \xymatrix{
\Spec{A}\ar@{_(->}[d]\ar[r] & \clK_{g,n}(\clG,\beta)\ar[d]\\
\Spec{B}\ar@{-->}[ur]\ar[r] & P_{g,n}.
}
\end{eqnarray}
Let us start by observing that the pullback  of $\clG$ to $\clC_B$ is a root gerbe by  Lemma \ref{Brauer_curve_over_artinian_ring}. Choose  a line bundle $L_B$ over $C_B$ and an isomorphism $\xi_B$ such that $\xi_B: f_B^*\clG\isomto \sqrt[r]{L_B/C_B}$. Let   $\iota: C_B|_A\hookrightarrow C_B$ denote the inclusion. There is a composite isomorphism 
\begin{eqnarray}\label{L_A_def_eq}
f_A^*\clG\isomto \varphi^* \iota^*f_B^*\clG \xrightarrow[\varphi^*\iota^*\xi_B]{\sim}
 \sqrt[r]{\varphi^*\iota^*L_B/C_A},
\end{eqnarray}
where the first arrow is canonical. Define $L_A:=\varphi^*\iota^* L_B$.  According to  this definition and  isomorphism (\ref{L_A_def_eq}), $\tilde{f}_A$ corresponds to a pair  $(\clN_A, \phi_A)$ where $\clN_A$ is a line bundle over $\clC_A$ and   $\phi_A$ is an isomorphism   such that $\phi_A: \clN_A^r\isomto L_A$. To give an extension of  $\tilde{f}_A$ to $B$ is equivalent to giving an extension $(\clN_B,\phi_B)$ of the pair $(\clN_A,\phi_A)$ and an isomorphism $\phi_B$ lifting $\phi_A$ such that $\phi_B: \clN_B^{\otimes r}\simeq L_B$.  By flatness we have an exact sequence
\begin{eqnarray}\label{sq-zero-ext-curves}
 1\to J\to \clO_{\clC_B} \to \clO_{\clC_A}\to 1,
\end{eqnarray}
where $J=I\otimes_B\clO_{\clC_B}$ is a square zero ideal in  $\clO_{\clC_B}$. By restricting to the subsheaves of invertible elements  and by taking the long exact cohomology sequence we get  
\begin{eqnarray}\label{ex-seq-Pic-curves}
  1\to H^1(\clC_A,\clO_{\clC_A})\otimes J\to \Pic{\clC_B}\to \Pic{\clC_A}\to 1,
\end{eqnarray}
where the last arrow is surjective because of dimensional reasons. Exactness on the left is due to the surjection $H^0(\clO_B^*)\to H^0(\clO_A^*)$. Choose a line bundle $\clN_B$ lifting $\clN_A$.  Let $\clS_B=\clN_B^{\otimes r}\otimes L_B^\vee$.  Then $\clS_B|_A\simeq \clO$, hence $\clS_B$ belongs to the subgroup $H^1(\clC_A,\clO_{\clC_A})\otimes J$. Since $J$ is a square zero ideal, it is not hard to see that $H^1(\clC_A,\clO_{\clC_A})\otimes J$ is divisible (e.g. by working with cocycles). We can therefore assume that $\clN_B^{\otimes r}\simeq L_B$. We need to show that any two line bundles $\clN'_B$, $\clN''_B$ lifting $\clN_A$ and satisfying the above condition are isomorphic. Indeed by assumption $\clN'_B\otimes{\clN''}_B^\vee$ is an $r$-torsion line bundle in $H^1(\clC_A,\clO_{\clC_A})\otimes J$. Again by using that $J$ is square zero one can see that this group does not contain torsion. Moreover any two choices of $(\clN_B,\phi_B)$ lifting $(\clN_A,\phi_A)$ and such that $\phi_B:\clN_B^{\otimes r}\isomto L_B$ define two isomorphic morphisms to $\sqrt[r]{L_B/C_B}$.
\end{pf}

Let $\clK:=\clK_{g,n}(\clG,\beta)^{\vec{g}}$ and $P:= P_{g,n}$. We observe that the relative inertia of the morphism $t:\clK_{g,n}(\clG,\beta)^{\vec{g}}\to P_{g,n}$, defined as $$I(\clK/P)=\clK \times_{\clK\times_P \clK}\clK,$$ contains as a subsheaf the \'etale sheaf $(\mu_r)_{\clK}$. Indeed, the automorphisms group  of an object  $[f:\clC_T\to\clG]$ over $T$ leaving $\clC_T$ fixed is $f^*(\mu_r)_\clG\simeq (\mu_r)_{\clC_T}$ bacause $\clG$ is a $\mu_r$-banded gerbe. Twisted curves are geometrically connected, therefore $\Gamma(\clC_T,\mu_r)=\Gamma(T,\mu_r)$. Whenever the the inertia of a finitely presented stack  $\clX$ contains an \'etale subgroup stack $G$ there exists a construction called {\em rigidification} \cite{AOV07}, which yields a canonical morphism $\clX\to \clX\thickslash G$ such that any morphism $f:\clX\to \clY$ whose relative inertia $I(\clX/\clY)$ contains $G$ factors through $\clX\to \clX\thickslash G$. \label{ridigif-page} Moreover $\clX \to \clX\thickslash G$  is an \'etale gerbe. In our case the rigidification of $\clK$ along $\mu_r$  is isomorphic to the relative coarse moduli space \cite{AOV08} for the morphism $t$. The relative coarse moduli space construction is recalled in Appendix \ref{morph-deg-appendix}.  By the above arguments we get the following

\begin{lem}\label{etale-fact-lem}
 The morphism $t:\clK_{g,n}(\clG,\beta)^{\vec{g}}\to P_{g,n}$ factors through
\begin{eqnarray}
 \clK_{g,n}(\clG,\beta)^{\vec{g}}\to \clK_{g,n}(\clG,\beta)^{\vec{g}}\thickslash \mu_r\to P_{g,n},
\end{eqnarray}
where the first arrow is a $\mu_r$-banded gerbe and the second arrow is representable and \'etale.
\end{lem}
\begin{pf}
 The second arrow is representable because  the rigidification coincides in this case with 
the  relative coarse moduli space.
The second arrow is also \'etale because the first arrow is \'etale and surjective, and being \'etale for a morphism
is a property \'etale local on the source.
\end{pf}


\section{Push-forward formula}
Consider again diagram (\ref{VFC_push_diag}). In this Section we prove the pushforward formula, Theorem \ref{push-thm}, which states that the pushforward along $p$ of $[\clK_{g,n}(\clG,\beta)^{\vec{g}}]^{vir}$ is a multiple of $[\barM_{g,n}(X,\beta)]^{vir}$ in $H_*(\barM_{g,n}(X,\beta),\bbQ)$.  We will show that the multiplicative factor is given by $r^{2g-1}$. The first step in the proof is to show that the natural perfect relative obstruction theory for $\clK_{g,n}(\clG,\beta)$ is quasi-isomorphic to the pull-back of the natural perfect
relative obstruction theory for $\barM_{g,n}(X,\beta)$. We refer the reader to Section \ref{vfc-subsection} for definitions and  notations.

\subsection{Comparison of Obstruction  Theories}

\begin{lem}\label{PROTHtwvsnontw}
Let $s: \clK_{g,n}(\clG,\beta)\to \barM_{g,n}(X,\beta)$ be the natural morphism.
 There is a natural isomorphism of objects in $\clD_{Coh}( \clK_{g,n}(\clG,\beta))$
\begin{equation*}
s^*E^\bullet \isomto \tilde{E}^\bullet.
\end{equation*}
\end{lem}

\begin{pf}
We will prove the statement for  $\tilde{E}^{\vee\bullet}=R\tilde{\pi}_*\tilde{f}^*T_\clG$ and $E^{\vee\bullet}=R\pi_*f^*T_X$. Consider the complex $Ls^*R\pi_*f^*T_X$ in $\clD_{coh}(\clK_{g,n}(\clG,\beta))$. It suffices to show that $Ls^*R\pi_*f^*T_X\simeq R\tilde{\pi}_*\tilde{f}^*T_\clG$. For this we consider the diagram 
\begin{equation*}
\xymatrix{
\clC\ar[dr]^\rho\ar[ddr]^{\tilde{\pi}}\ar[drr]\ar[rr]^{\tilde{f}} & & \clG\ar[dr]^\epsilon\\
&p^*C\ar[d]^{\pi'}\ar[r]^{s'} &C\ar[d]^\pi\ar[r]^f & X\\
&\clK_{g,n}(\clG, \beta)\ar[r]^s & \overline{M}_{g,n}(X, \beta).
}
\end{equation*}
Observe that $\epsilon^*T_X\simeq T_\clG$. Also we have $R\rho_*L\rho^*\simeq Id$ because the map $\rho$ is the relative coarse moduli space for the map $\tilde{\pi}$. The  arrow $s$ is flat, because it is the composition of a flat morphism with an \'etale morphism. The arrow $\pi$ is flat because it  is the structure morphism of the universal curve. 
 Moreover  the square in the  above diagram is cartesian, hence we calculate (using \cite[Proposition 13.1.9]{LMBca})
\begin{equation*}
\begin{split}
Ls^*R\pi_*f^*T_X& \simeq R\pi'_{*}Lp'^*f^*T_X\\
&\simeq R\pi'_{*}R\rho_*L\rho^*Lp'^*f^*T_X\\
&\simeq R\pi'{*}R\rho_*\tilde{f}^*\epsilon^*T_X\\
&\simeq R\pi'_{*}R\rho_*\tilde{f}^*T_\clG\\
&\simeq R\tilde{\pi}_*\tilde{f}^*T_\clG.
\end{split}
\end{equation*}
Since $s^*$  is  exact,  we write $s^*$ for $Ls^*$. 
\end{pf}

\begin{lem}
 The diagram of morphisms in $\clD_{coh}(\clK_{g,n}(\clG,\beta))$
\begin{eqnarray}
 \xymatrix{
s^*E^\bullet \ar[r]\ar[d] & s^*L_{\barM_{g,n}(X,\beta)/\MM_{g,n}}^\bullet\ar[d]\\
\tilde{E}^\bullet \ar[r] & L_{\clK_{g,n}(\clG,\beta)/\MM_{g,n}^{tw}}^\bullet
}
\end{eqnarray}
is commutative. 
\end{lem}
\begin{pf}
The obstruction theories above are determined  by adjunction from the morphisms obtained by composing the arrows in the following diagram
\begin{eqnarray}\label{obstr-thy-comm-diag}
 \xymatrix{
s^*f^*L^\bullet_X \ar[r]\ar[d] & \tilde{f}^* L_\clG\ar[d]\\
s^*L_{C}\ar[r] \ar[d] & L^\bullet_{\clC}\ar[d]\\
s^*L^\bullet_{C/C'}\ar[r] & L^\bullet_{\clC/\clC'}\\
s^*\pi^*L^\bullet L_{\barM_{g,n}(X,\beta)/\MM_{g,n}} \ar[u]\ar[r] & \tilde{\pi}^* L_{\clK_{g,n}(\clG,\beta)/\MM_{g,n}^{tw}}^\bullet\ar[u]
}
\end{eqnarray}
Each arrow in diagram (\ref{obstr-thy-comm-diag}) is commutative since every arrow is part of a transitivity exact sequence involving cotangent complexes. Commutativity follows from functorial properties of transitivity exact sequences  (cfr. \cite[2.1.5]{Ill_I}). 
\end{pf}

\begin{cor}\label{obs-thy-pullback-coro}
The induced morphism $s^*E^\bullet\to L_{\clK_{g,n}(\clG,\beta)/\MM_{g,n}^{tw}}$ coincides with the orbifold GW obstruction theory
$\tilde{E}^\bullet\to  L_{\clK_{g,n}(\clG,\beta)/\MM_{g,n}^{tw}}$.
\end{cor}
\begin{cor}\label{obs-thy-pullback-coro-comp}
For any $\beta$-admissibe vector $\vec{g}$, the induced morphism $s^*E^\bullet\to L_{\clK_{g,n}(\clG,\beta)^{\vec{g}}/\MM_{g,n}^{tw}}$ coincides with the orbifold GW obstruction theory $\tilde{E}^\bullet\to  L_{\clK_{g,n}(\clG,\beta)^{\vec{g}}/\MM_{g,n}^{tw}}$.
\end{cor}
\begin{pf}
$ \clK_{g,n}(\clG,\beta)^{\vec{g}}\hookrightarrow \clK_{g,n}(\clG,\beta)$ is an open immersion.
\end{pf}

\subsection{Proof of the pushforward  formula}\label{proof-push-sect}
By Corollary \ref{obs-thy-pullback-coro-comp}  the natural perfect relative obstruction theory $\tilde{E}^\bullet$  over $\clK_{g,n}(\clG,\beta)^{\vec{g}}$ is quasi-isomorphic to the  pull-back of the natural perfect relative obstruction theory $E^\bullet$ over $\barM_{g,n}(X,\beta)$. We know that such obstruction theories admit global resolutions. We choose a length 2 complex of vector bundles  $F^{\bullet}$ quasi isomorphic to $E^\bullet$. Its pullback to $\clK_{g,n}(\clG,\beta)^{\vec{g}}$, that we denote by $\tilde{F}^\bullet$, is a global resolution of $\tilde{E}^\bullet$. We denote the dual complexes by $\tilde{F}_\bullet$ and $F_\bullet$. We have a cartesian diagram (put $\clK=\clK_{g,n}(\clG,\beta)^{\vec{g}}$, $\barM=\barM_{g,n}(X,\beta)$)
\begin{eqnarray}\label{cones-cart-diag}
 \xymatrix{
\tilde{C}\ar@{^{(}->}[d]\ar[r]^{p''}\ar@{}[rd]|{\square} & C\ar@{^{(}->}[d]\\
\tilde{F}_1\ar[d]_{\tilde{\pi}}\ar@{}[rd]|{\square}\ar[r]_{p'} & F_1\ar[d]^{\pi}\\
\clK\ar[r]_{p} & \barM
}
\end{eqnarray}
where $\tilde{C}= \tilde{F}_1\times_{[\tilde{F}_1/[\tilde{F}_0]}\mathfrak{C}_{\clK/\MM^{tw}}$
and $C=F_1\times_{[F_1/F_0]}\mathfrak{C}_{\barM/\MM}$. We observe that $\tilde{F}_1$ and $F_1$ are vector bundles. Hence the intersections with the zero sections provide isomorphisms  $0_{\clK}^!: A_*(\tilde{F}_1)_\bbQ\isomto A_*(\clK)_\bbQ$ and  $0_{\barM}^!: A_*(F_1)_\bbQ\isomto A_*(\barM)_\bbQ$.  By compatibility of Gysin morphisms with proper pushforward we get that for any
$\alpha\in A_*(\tilde{F}_1)_\bbQ$
\begin{eqnarray}
 p_* \ 0_{\clK}^!\ \alpha= 0_{\barM}^!\  p'_*\ \alpha.
\end{eqnarray}
We argue in the spirit of  \cite{Cost03} that $p_*[\clK]^{vir}=d[\barM]^{vir}$ for some $d\in\bbQ$ if and only if 
\begin{eqnarray}
p'_*[\tilde{C}]=d[C] \in A_*(F_1)_\bbQ.
\end{eqnarray}
We will show that $\tilde{C}$ is the pullback of $C$ and that
\begin{eqnarray}\label{push-cones-eq}
p'_*[\tilde{C}]=p'_*p'^*[C]=r^{2g-1} [C] \in A_*(F_1)_\bbQ.
\end{eqnarray}

\begin{prop}
Consider diagram (\ref{VFC_push_diag}). The intrinsic normal cone $\mathfrak{C}_{\clK/\MM^{tw}}$ is isomorphic to the  pullback of the relative instrinsic normal cone $\mathfrak{C}_{\barM/\MM}$, where $\clK=\clK_{g,n}(\clG,\beta)^{\vec{g}}$, $\barM=\barM_{g,n}(X,\beta)$, $\MM^{tw}=\MM^{tw}_{g,n}$ and $\MM=\MM_{g,n}$.
\end{prop}
\begin{pf}
Note first that the relative instrinsic normal cone $\mathfrak{C}_{P_{g,n}/\MM^{tw}}$ is isomorphic to the pullback of $\mathfrak{C}_{\barM/\MM}$ because $q$ is flat. This follows from \cite[Proposition 7.2 ]{BFinc}. We claim that $\mathfrak{C}_{\clK/\MM^{tw}}$ is isomorphic to $t^*\mathfrak{C}_{P_{g,n}/\MM^{tw}}$.  By the functorial properties  of the intrinsic normal cone there is a natural morphism
\begin{eqnarray}\label{pull-back-nc-etale}
 \mathfrak{C}_{\clK/\MM^{tw}}\to t^*\mathfrak{C}_{P_{g,n}/\MM^{tw}}.
\end{eqnarray}
Being an isomorphism is a local property. We will show that such a morphism is an isomorphism in local charts. Any local embedding of $\clK$ over $\MM^{tw}_{g,n}$ is also a local embedding for $P_{g,n}$ over $\MM^{tw}_{g,n}$ because $t$ is \'etale. Let $(U,M)$ such a local embedding. Restrict the natural morphism (\ref{pull-back-nc-etale}) to $U$. By the local description of the intrinsic relative normal cone we get
\begin{eqnarray}
\xymatrix{
 \mathfrak{C}_{\clK/\MM^{tw}}|_U\ar[r]\ar[d]^{\sim} & t^*\mathfrak{C}_{P_{g,n}/\MM^{tw}}|_U\ar[d]^{\sim}\\
[C_U/T_{M/\MM^ {tw}}]\ar[r]^{\sim} &  [C_U/T_{M/\MM^ {tw}}].
}
\end{eqnarray}
\end{pf}

In order to prove the equality  in (\ref{push-cones-eq}) we will show that $p$, hence $p'$ is of pure degree. We will use the facts recalled in Appendix \ref{morph-deg-appendix}.
We make a few preliminary remarks. The morphism $t:\clK_{g,n}(\clG,\beta)^{\vec{g}} \to P_{g,n}^{\vec{g}}$ is \'etale by Proposition \ref{t-etaleness-prop}. The morphism
$g:\MM_{g,n}^{tw}\to\MM_{g,n}$ is flat by \cite{OLogCurv}. Moreover $p:\clK_{g,n}(\clG,\beta)^{\vec{g}}\to \barM_{g,n}(X,\beta) $ is proper  because it is a morphism between proper stacks and is quasi finite because e.g. it is a composition of two quasi finite morphisms. Indeed, the image of $\clK_{g,n}(\clG,\beta)^{\vec{g}}$ in $\MM_{g,n}^{tw}$ factors through the open substack where the orders of isotropy groups of points in the boundary divisors divide $r$. Such an open substack is of finite type and quasi-finite over $\MM_{g,n}$. Therefore $q'$ is also quasi-finite. The morphism $t$ is quasi-finite because it is \'etale. Finally, $p=q'\circ t$, hence $p$ is quasi-finite.

\begin{prop}\label{comp-deg-prop}
With the notations in diagram (\ref{cones-cart-diag})
\begin{eqnarray}\label{push-con-eq-prop}
 p'_*[\tilde{C}]=r^{2g-1}[C].
\end{eqnarray}
\end{prop}

\begin{pf}
Let us denote by $C^r$ the reduced substack of an irreducible component of $C$. Let $\tilde{C}'=p^{-1}(C^r)$ and let $\tilde{C}^r$ be the reduced substack of $\tilde{C}$ and $\tilde{C}'$. Let $m(C)$, $m(\tilde{C}')$ denote the geometric multiplicities. Equation (\ref{push-con-eq-prop}) can be written as
\begin{eqnarray}
 p'_*p'^*[C] = m(C)m(\tilde{C}')p'_*[\tilde{C}^r]=m(C)m(\tilde{C}')\mbox{deg}(\tilde{C}^r/C^r)[C^r]\nonumber \\ = \mbox{deg}(\tilde{C}'/C^r)[C].
\end{eqnarray}
We compute the degree of $\tilde{C}'$ over $C^r$ by using the characterization given in Appendix \ref{morph-deg-appendix}. Assume that $C^r$ factors through an integral substack $\barM(\tau)^{r}$ of $\barM_{g,n}(X,\beta)$ such that  its generic geometric points  correspond to stable maps with domain curve of dual graph $\tau$. Let $\clK(\tau)$ be the preimage of $\barM(\tau)^{r}$ along $p$. The dimension of the ring $K(\clK)$ over $K(\barM(\tau)^{r})$,  $[K(\clK):K(\barM(\tau)^{r})]$,  is equal to  the number of points weighted with multiplicity in a generic geometric fiber of $\overline{p}:\overline{\clK}(\tau)\to \barM(\tau)^{red}$, which is the relative coarse moduli space of $p:\clK(\tau)\to \barM(\tau)^{red}$ (see Appendix \ref{morph-deg-appendix}). 
Let $E_\tau^{loop}\subset E_\tau$ be the subset of looping edges. Let $P(\tau)=P_{g,n}\times_{\barM_{g,n}(X,\beta)}\barM(\tau)^{red}$. Let $\overline{x}: \Spec{\bbC} \to \barM(\tau)^{red}$ be a generic  geometric point. Let $P(\tau)(\overline{x})$ be the fiber of $P(\tau)\to \barM(\tau)^{red}$ over $\overline{x}$. Then
\begin{eqnarray}
 \clK(\tau)\times_{\barM(\tau)^{red}}\Spec{\bbC}\simeq  \clK(\tau)\times_{P(\tau)} P(\tau)(\overline{x}).
\end{eqnarray}
Let $A(\tau)$ be the set of all $\vec{g}$-compatible gerby graphs over $\tau$. For any $\ttau$ in $A(\tau)$ and for any $i=1,..,\abs{E_\ttau}$ let $s_i(\ttau)=\gamma(e_i)$, $e_i\in E_\ttau$. Then
\begin{eqnarray}
 P(\tau)(\overline{x})= \coprod_{\ttau\in A(\tau)} [\Spec{\bbC[x_1,...,x_r]}/(x_i^{s_i(\ttau)})/\mu(\ttau)]:=\coprod_{\ttau\in A(\tau)} W(\ttau),
\end{eqnarray}
where $\mu(\ttau)=\mu_{s_1(\ttau)}\times ...\times \mu_{s_r(\ttau)}$  and $(\epsilon_1,..,\epsilon_r)$ in $\mu(\ttau)$ acts on $(x_1,..,x_r)$ taking it to $(\epsilon_1 x_1,..,\epsilon_r x_r)$. 
 Let  $J(\tau)\subset A(\tau)$ be the set  of $\vec{g}$-compatible gerby graphs $\ttau$ over $\tau$ such that for any $e_i\in E_\ttau$, $\gamma(e_i)$ divides $r$. Let $I(\ttau)$ be the number of non-isomorphic line bundles over a twisted curve of gerby dual graph $\ttau$. %
For any $\ttau\in J(\ttau)$  there is a cartesian diagram
\begin{eqnarray}
 \xymatrix{
\coprod_{l(\ttau)} (B\mu_r)_{l(\ttau)}\ar[r] \ar[dd] & \coprod_{l(\ttau)} Z(\ttau)_{l(\ttau)} \ar[d]\ar[r] & \clK(\tau)(\overline{x})\ar[d]\\
& \coprod_{l(\ttau)} W(\ttau)_{l(\ttau)}\ar[d]\ar[r] & \clK(\tau)\thickslash\mu_r (\overline{x})\ar[d]\\
\Spec{\bbC} \ar[r] & W(\ttau)\ar[r] & P(\tau)(\overline{x})
}
\end{eqnarray}
where $1\leq l(\ttau)\leq I(\ttau)$, $(B\mu_r)_{l(\ttau)}=(B\mu_r)$, $Z(\ttau)_{l(\ttau)}=Z(\ttau)$ and $W(\ttau)_{l(\ttau)}=W(\ttau)$. By Lemma \ref{etale-fact-lem} the upper vertical arrows are the structure morphisms of a $\mu_r$-gerbe and the lower vertical arrows are \'etale and representable. In our case, as can be checked by using groupoid presentations\footnote{Artinian local rings with residue field $\bbC$ admit only trivial \'etale covers}, the lower central arrow is a trivial cover.
The upper left corner is a disjoint union of $I(\ttau)$ trivial gerbes. We recall that $I(\ttau)$ is the number of non isomorphic $r$-torsion line bundles over a twisted curve of dual graph $\ttau$ associated to a fixed $\vec{g}$ and with fibers carrying faithful representations of the isotropy groups of the nodes. Connected components of $\clK(\tau)(\overline{x})\times_{P(\tau)(\overline{x})} W(\tilde{\tau})$ are labeled by the corresponding connected component in $\coprod_{l(\ttau)} (B\mu_r)_{l(\ttau)}$. Note that for any $\ttau$ the (relative) coarse moduli space of $W(\ttau)$ (over $\Spec{\bbC}$) is $\Spec{\bbC}$ and for any $\ttau$ and any $l(\ttau)$ the same is true for $Z(\ttau)_{l(\ttau)}$. Therefore the  following  diagram is cartesian
 \begin{eqnarray}\label{rel-cms-diag}
 \xymatrix{
\coprod_{\ttau\in J(\tau)}\coprod_{l(\ttau)} Z(\ttau)_{l(\ttau)} 
\ar[d]_{\pi_{\overline{x}}} \ar@{}[rd]|{\square}\ar[r]& \clK(\tau)\ar[d]^{\pi}\\
 \coprod_{\ttau\in J(\tau)}\coprod_{l(\ttau)} (\Spec{\bbC})_{l(\ttau)} \ar[d]_{\overline{p}_{\overline{x}}} \ar[r]\ar@{}[rd]|{\square} &  \overline{\clK}(\tau)\ar[d]^{\overline{p}}\\
  \Spec{\bbC}\ar[r]_{\overline{x}} & \barM(\tau)^{r}
 }
\end{eqnarray}
Here $\overline{p}:\overline{\clK}\to \barM(\tau)^{r}$ is the relative coarse moduli space for $p$ and  $\overline{p}_x$ is the relative coarse moduli space of  $\overline{p}_x \circ \pi_{\overline{x}}$ because the formation of the relative coarse moduli space commutes with arbitrary base change for tame stacks. The number of points in the fiber of $\overline{p}_{\overline{x}}$ is computed in Lemma \ref{number-points-fiber-lemm} below. It is equal to $r^{2g}$ and  it does not depend on $\tau$. According to the degree formula in Appendix \ref{morph-deg-appendix}, in order to compute the degree of $p$ we need to multiply by the number $\delta(\clK_{g,n}(\clG,\beta))/\delta(\barM_{g,n}(X,\beta))$ in case there are non trivial generic stabilizers. Note that for any $T\to \clK_{g,n}(\clG,\beta)$ corresponding to an object $[\clC\to \clG]$ over $[C\to X]$, there is a surjection of sheaves of groups  $\underline{Aut}_T(\clC\to \clG)\to\underline{Aut}_T(C\to X)$. Hence  $\delta(\clK_{g,n}(\clG,\beta))/\delta(\barM_{g,n}(X,\beta))$ is equal to the degree of the relative inertia $I(\clK/\barM)$. This  is in turn equal to $\mbox{deg}\,I(\clK/\overline{\clK})=r$, because $\clK\to \overline{\clK}$ is a $\mu_r$-gerbe. By putting all together we get
\begin{eqnarray}
 \mbox{deg } p= r^{2g-1}.
\end{eqnarray}
\end{pf}

\begin{lem}\label{number-points-fiber-lemm}
Consider the morphism $\overline{p}:\overline{\clK}(\tau)\to \barM(\tau)$ defined in Proposition \ref{comp-deg-prop}. Let $\overline{x}:\Spec{\bbC}\to \barM(\tau)^{r}$ be a generic geometric point. Then the number of geometric points of $\overline{\clK}(\tau)$ over  $\overline{x}$ is equal to $r^{2g}$.
\end{lem}

\begin{pf}
We first observe that geometric points of $\overline{\clK}(\tau)$ are in bijection with geometric points of $\clK(\tau)$. Therefore we have to count the number of (non isomorphic) twisted curves $\clC$ over $C$ with $\vec{g}$-compatible gerby dual graph  over $\tau$  admitting twisted stable maps to $\clG$. Moreover, for any  such $\clC$ we have to count the number of non-isomorphic twisted stable maps to $\clG$. Note first that a twisted curve over $\Spec{\bbC}$ is determined by its coarse moduli space and by the order of the isotropy groups of its special points\footnote{This can be understood e.g.  by considering the characterization of twisted curves in terms of log-twisted curves given in \cite{OLogCurv} and by observing that the only isomorphism class of locally free log structure over $\Spec{\bbC}$ corresponds to a monoid of the form $\bbN^{r}$, where $r\in\bbN$ is equal  to the number of nodes of  the coarse moduli space $C$}. Let us start by computing  the number of $\vec{g}$-compatible gerby graphs $\ttau$ over $\tau$ such that the corresponding $\clC$ admits a twisted stable map to $\clG$. By Proposition \ref{geom-lift-prop}
we know that any twisted curve with $\vec{g}$-compatible $\ttau$ such that for any $e_i\in E_\ttau^{loop}$ $\gamma(e_i)$ divides $r$ admits twisted stable maps to $\clG$. Choose a line bundle $L$ such that $f^*\clG\simeq \sqrt[r]{L/C}$. We fix a set of line bundles satisfying the properties listed in Property \ref{canonical-tw-lb-nodes-property} and we use the same notations we used there. After making this choice, the admissible vector $\vec{g}$ determines a line bundle $\clS$ over $\clC$ such that $\mbox{deg}\, L\otimes \clS^{-r}$ is multiple of $r$. Explicitly, $\clS\simeq \otimes_{i=1}^n \clT_i^{m_i}\bigotimes \otimes^{n.l.}_{j\in E_\ttau} \clQ_j^{t_j}$ where for all $e_i\in F_\ttau\setminus E_\ttau$,  $(m_i,\gamma(e_i))=1$ and for all $e_j\in E_\ttau\setminus E_\ttau^{loop}$, $(\gamma(e_j),t_j)=1$. Any line bundle defined as
\begin{eqnarray}\label{possible-torsion-choices}
\clN_{\vec{s}}= \clS\otimes\clR_{\vec{s}}=\clS \otimes \otimes^{e_l\in E_\ttau^{loop}} \clQ_l^{s_l},\quad (s_l,\gamma(e_l))=1\ \ \forall e_l\in E_\ttau^{loop}
\end{eqnarray}
corresponds to a twisted stable map $\clC\to\sqrt[r]{L/C}\simeq f^*\clG$. Indeed by definition $\clR_{\vec{s}}^{\otimes r}\simeq \clO_\clC$. The condition $(s_l,\gamma(e_l))=1$ as usual implies the representability of the corresponding morphism. By tensoring $\clS$ with different possible choices of line bundles of the form $\clR_{\vec{s}}$  we get different isomorphism classes of twisted stable maps. The number of choices of $\clR_{\vec{s}}$ is computed by
\begin{eqnarray} 
N(r_1,..,r_k)=\prod_{l=1}^k \phi(r_l),
\end{eqnarray}
where $k=\abs{E_\ttau^{loop}}$, $r_l=\gamma(e_l)$ for any $e_l\in E_\ttau^{loop}$ and  $\phi(-)$ is the Euler totient function\footnote{For any $n\in\bbN$, the Euler totient function $\phi(n)$ is defined to be the number of integers less than $n$ and coprime with $n$. The following property holds: $\sum_{d|n}\phi(d)=n$.}. Moreover, by tensoring any $\clN_{\vec{s}}$ with different possible choices of   $r$-torsion line bundles in  $\Pic{C}$ we get $r^{2g-b_1(\tau)}$  non isomorphic twisted stable maps. Finally, in order to get the total number of geometric points over $\overline{x}$ we have to sum over all $\ttau\in J(\ttau)$ with $J(\ttau)$ defined in Proposition \ref{comp-deg-prop}. This is the same as summing the numbers $N(r_1,..,r_n)$ over all possible $n$-tuples $(r_1,..,r_n)$ of integers dividing $r$. Let $\abs{E_\tau^{loop}}=k$. By the properties of Euler totient function we get
\begin{eqnarray}
\tilde{N}= \sum_{r_1,..,r_k|r} N(r_1,..,r_k)=\sum_{r_1,..,r_k|r} \prod_{l=1}^k \phi(r_l)=r^k=r^{b_1(\tau)}.
\end{eqnarray}
Finally, we get the total number of points over $\overline{x}$ as
\begin{eqnarray}
 N_{tot}=\tilde{N}\times r^{2g-b_1(\tau)}=r^{2g}.
\end{eqnarray}
Note that this number does not depend on the dual graph $\tau$.
\end{pf}

Combining the above discussions, we get the desired pushforward formula.
\begin{thm}\label{push-thm}
Let $\vec{g}$ be a $\beta$-admissible vector. Let $p: \clK_{g,n}(\clG, \beta)^{\vec{g}}\to\barM_{g,n}(X,\beta)$ be the  morphism defined  in diagram (\ref{VFC_push_diag}). Then
\begin{eqnarray}
 p_*[\clK_{g,n}(\clG, \beta)^{\vec{g}}] ^{vir}=r^{2g-1}[\barM_{g,n}(X,\beta)]^{vir} \in H^*(\barM_{g,n}(X,\beta),\bbQ).
\end{eqnarray}
\end{thm}
\begin{pf} 
It follows from arguments at the beginning of Section \ref{proof-push-sect} and from Proposition
\ref{comp-deg-prop}.
 \end{pf}
 
\section{Orbifold Gromov-Witten theory of banded gerbes}
In this Section we examine the Gromov-Witten invariants of the gerbe $\clG$ using results in previous sections. In particular in Theorem \ref{decomp_any_genus} we prove the decomposition conjecture for $\clG$. 

\subsection{Orbifold Gromov-Witten invariants}
Let $$\epsilon: \clG \to X$$  be a $G$-banded gerbe with $G$ a finite abelian group over $X$. Let
$\alpha\in H^2(X,G)$ be the isomorphism class of $\clG$ (as a $G$-banded gerbe). Since  the 
gerbe is $G$-banded   there is a canonical isomorphism
\begin{eqnarray}
 I\clG=G\times_X \clG\simeq \coprod_{g\in G}\clG_{g},
\end{eqnarray}
where $\clG_{g}$ is a  root gerbe isomorphic to  $\clG$. Let $\epsilon_{g}:\clG_{g}\to X$
be the induced morphism.  On each component there is an isomorphism between the rational cohomology groups 
$$\epsilon_{g}^{*}: H^{*}(X,\mathbb{Q})\overset{\simeq}{\longrightarrow} H^{*}(\clG_{g}, \mathbb{Q}).$$

To simplify notation, from now on we assume $G=\mu_r$ for some $r\in \bbN$. The discussion for the general $G$ requires only notational changes. The genus $0$ case is spelled out explicitly in \cite[Appendix A]{AJT09-gerbes}. The higher genus case is similar.

Let $\vec{g}=(g_1,...,g_n)$ be a $\beta$-admissible vector. There are evaluation maps $$ev_i: \mathcal{K}_{g,n}(\clG,\beta)^{\vec{g}}\to \bar{I}(\clG)_{g_i},$$
where $\bar{I}(\clG)_{g_i}$ is a component of  the  {\em rigidified inertia stack} $\bar{I}(\clG)=\cup_{g\in \mu_r} \bar{I}(\clG)_g$. Although the evaluation maps  $ev_i$ do not take values in  $I\clG$, as explained in \cite{AGV06}, Section 6.1.3, one can still define a pull-back map at cohomology level,
$$ev_i^*: H^*(\clG_{g_i}, \mathbb{Q})\to H^*(\mathcal{K}_{g,n}(\clG,\beta)^{\vec{g}}, \mathbb{Q}).$$
Given $\delta_i\in H^*(\clG_{g_{i}},\mathbb{Q})$ for $1\leq i\leq n$ and   integers $k_i\geq 0,1\leq i\leq n$,  one can define descendant orbifold Gromov-Witten invariants 
$$\langle \delta_1\bar{\psi}_1^{k_1},\cdots, \delta_n\bar{\psi}_n^{k_n}\rangle_{g, n, \beta}^\clG:=\int_{[\mathcal{K}_{g,n}(\clG,\beta)^{\vec{g}}]^{vir}}\prod_{i=1}^n ev_i^*(\delta_i)\bar{\psi}_i^{k_i},$$
where $\overline{\psi}_i$ are the pullback of the first  Chern classes  of the tautological line bundles  over $\overline{M}_{g,n}(X,\beta)$ (which by abuse of notation we also denote by $\bar{\psi}_i$). 

For  classes $\delta_i\in H^*(\clG_{g_i},\mathbb{Q})$, 
set $\overline{\delta}_i=(\epsilon_{g_{i}}^{*})^{-1}(\delta_i)$. Descendant 
Gromov-Witten invariants  $\langle\overline{ \delta}_1\bar{\psi}_1^{k_1},\cdots,\overline{\delta}_n\bar{\psi}_n^{k_n}\rangle_{0, n, \beta}^X$ of $X$ are similarly defined. 
Theorem \ref{push-thm} implies the following comparison result.
\begin{thm}\label{GW-inv1}
\begin{equation*}
 \langle \delta_1\bar{\psi}_1^{k_1},..., \delta_n\bar{\psi}_n^{k_n}\rangle_{g,n,\beta}^\clG=r^{2g-1} \langle \overline{\delta}_1\bar{\psi}_1^{k_1},\cdots,\overline{\delta}_n\bar{\psi}_n^{k_n}\rangle_{g,n,\beta}^X.
\end{equation*}
Moreover, if $\vec{g}$ is not admissible, then the Gromov-Witten invariants of $\clG$ vanish.
\end{thm}

\begin{pf}
Denote by $\overline{ev}_i: \overline{M}_{g,n}(X, \beta)\to X$ the $i$-th evaluation map. Using the definition of $ev_i^*$ one can check that $ev_i^*(\delta_i)=p^*\overline{ev}_i^*(\overline{\delta}_i)$. Note also that $p^*\bar{\psi}_i=\bar{\psi}_i$. Thus using Theorem \ref{push-thm} we have 
\begin{equation*}
\begin{split}
\langle \delta_1\bar{\psi}_1^{k_1},..., \delta_n\bar{\psi}_n^{k_n}\rangle_{g,n,\beta}^\clG&=\int_{[\mathcal{K}_{g,n}(\clG,\beta)^{\vec{g}}]^{vir}}\prod_{i=1}^n ev_i^*(\delta_i)\bar{\psi}_i^{k_i}\\ 
&=\int_{[\mathcal{K}_{g,n}(\clG,\beta)^{\vec{g}}]^{vir}}\prod_{i=1}^n p^*\overline{ev}_i^*(\overline{\delta}_i)\bar{\psi}_i^{k_i}\\
&=\int_{[\mathcal{K}_{g,n}(\clG,\beta)^{\vec{g}}]^{vir}}\prod_{i=1}^n p^*(\overline{ev}_i^*(\overline{\delta}_i)\bar{\psi}_i^{k_i})\\
&=r^{2g-1}\int_{[\overline{M}_{g,n}(X,\beta)]^{vir}}\prod_{i=1}^n \overline{ev}_i^*(\overline{\delta}_i)\bar{\psi}_i^{k_i}\\
 &=r^{2g-1} \langle \overline{\delta}_1\bar{\psi}_1^{k_1},\cdots,\overline{\delta}_n\bar{\psi}_n^{k_n}\rangle_{g,n,\beta}^X.
\end{split}
\end{equation*}
\end{pf}
\subsection{Decomposition of Gromov Witten theory}
In the following we use complex numbers $\bbC$ as coefficients for the cohomology. For $\overline{\alpha}\in H^*(X, \mathbb{C})$ and an irreducible representation $\rho$ of $G$, we define $$\overline{\alpha}_\rho:=\frac{1}{r}\sum_{g\in G} \chi_{\rho}(g^{-1})\epsilon_g^*(\overline{\alpha}),$$ where $\chi_\rho$ is the character of $\rho$. The map $(\overline{\alpha},\rho)\mapsto \overline{\alpha}_\rho$  clearly defines an additive isomorphism
\begin{equation}\label{char_table_isom}
\bigoplus_{[\rho]\in \widehat{G}}H^*(X)_{[\rho]}\simeq H^*(I\clG, \mathbb{C}),
\end{equation}
where $\widehat{G}$ is the set of isomorphism classes of irreducible representations of $G$, and for $[\rho]\in\widehat{G}$ we define $H^*(X)_{[\rho]}:= H^*(X, \mathbb{C})$.

Theorem \ref{GW-inv1} together with orthogonality relations of characters of  $G$ implies the following 

\begin{thm}\label{GW-inv2}
Given $\overline{\alpha}_1, ...,\overline{\alpha}_n\in H^*(X, \bbQ)$ and integers $k_1,..., k_n\geq 0$, we have 
\begin{equation*}
\begin{split}
&\langle \overline{\alpha}_{1\rho_1}\bar{\psi}_1^{k_1},..., \overline{\alpha}_{n\rho_n}\bar{\psi}_n^{k_n} \rangle_{g,n,\beta}^{\clG}\\
=&\begin{cases}
r^{2g-2}\langle \overline{\alpha}_1\bar{\psi}_1^{k_1},\cdots,\overline{\alpha}_n\bar{\psi}_n^{k_n}\rangle_{g,n,\beta}^X\chi_\rho(\exp(\frac{-2\pi \sqrt{-1}(\alpha\cap\beta)}{r})) &\text{if }\rho_1=\rho_2=...=\rho_n=:\rho,\\
0&\text{otherwise}\,.
\end{cases}
\end{split}
\end{equation*}
(Recall that $\alpha\in H^2(X, \mu_r)$ is the class of the gerbe $\clG\to X$.)
\end{thm}

\begin{pf}
By our definition we have 
\begin{equation*}
\begin{split}
&\langle \overline{\alpha}_{1\rho_1}\bar{\psi}_1^{k_1},..., \overline{\alpha}_{n\rho_n}\bar{\psi}_n^{k_n} \rangle_{g,n,\beta}^{\clG}\\
=&r^{2g-1}\sum_{g_1,...,g_n\in \mu_r}\prod_{i=1}^n\chi_{\rho_i}(g_i^{-1})\langle \prod_{i=1}^n \epsilon_{g_i}^*(\overline{\alpha}_i)\bar{\psi}_i^{k_i} \rangle_{g,n,\beta}^\clG.
\end{split}
\end{equation*}
The term associated to $\vec{g}:=(g_1,...,g_n)$ in the above sum vanishes unless $\vec{g}$ is an admissible vector. This implies that $\prod_{i=1}^n g_i=\exp(\frac{2\pi \sqrt{-1} (\alpha\cap \beta)}{r})$. We rewrite this as $g_n^{-1}=\exp(\frac{-2\pi\sqrt{-1}(\alpha\cap\beta)}{r} )\prod_{i=1}^{n-1} g_i$. Substitute this into above equation and use Theorem \ref{GW-inv1} to get 
\begin{equation*}
\begin{split}
&\langle \overline{\alpha}_{1\rho_1}\bar{\psi}_1^{k_1},..., \overline{\alpha}_{n\rho_n}\bar{\psi}_n^{k_n} \rangle_{g,n,\beta}^{\clG}\\
=&r^{2g-1}\sum_{g_1,...,g_{n-1} \in \mu_r}\chi_{\rho_n}(\exp(\frac{-2\pi \sqrt{-1}(\alpha\cap\beta)}{r}))\left(\prod_{i=1}^{n-1}\chi_{\rho_i}(g_i^{-1})\chi_{\rho_n}(g_i)\right)\frac{1}{r}\langle \prod_{i=1}^n \overline{\alpha}_i\bar{\psi}_i^{k_i} \rangle_{g,n,\beta}^X.
\end{split}
\end{equation*}
Applying the orthogonality condition 
$$\frac{1}{r}\sum_{g\in \mu_r}\chi_\rho(g^{-1})\chi_{\rho'}(g)=\delta_{\rho, \rho'},$$
we find 
\begin{equation*}
\begin{split}
&\langle \overline{\alpha}_{1\rho_1}\bar{\psi}_1^{k_1},..., \overline{\alpha}_{n\rho_n}\bar{\psi}_n^{k_n} \rangle_{g,n,\beta}^{\clG}\\
=&\frac{1}{r}\chi_{\rho_n}(\exp(\frac{-2\pi \sqrt{-1}(\alpha\cap\beta)}{r}))\left(\prod_{i=1}^{n-1}\delta_{\rho_i, \rho_n}\right) r^{2g-1}\langle \prod_{i=1}^n \overline{\alpha}_i\bar{\psi}_i^{k_i} \rangle_{g,n,\beta}^X.
\end{split}
\end{equation*}
The result follows.
\end{pf}

We now reformulate this in terms of generating functions. Let $$\{\overline{\phi}_i\, |\, 1\leq i\leq \text{rank}H^*(X, \mathbb{C})\}\subset H^*(X, \mathbb{C})$$ be an additive basis. According to the discussion above, the set $$\{\overline{\phi}_{i\rho}\,|\, 1\leq i\leq \text{rank}H^*(X, \mathbb{C})\}, [\rho]\in \widehat{G}$$ is an additive basis of $H^*(I\clG, \mathbb{C})$. Recall that the genus $g$ descendant potential of $\clG$ is defined to be 
\begin{eqnarray}
&\mathcal{F}^g_{\clG}(\{t_{i\rho, j}\}_{1\leq i\leq \text{rank}H^{*}(X,\mathbb{C}), \rho\in \widehat{G}, j\geq 0}; Q):= & \nonumber \\
& \sum_{\overset{n\geq 0, \beta\in H_2(X,\mathbb{Z})}{i_1,...,i_n; \rho_1,...,\rho_n; j_1,...,j_n}}\frac{Q^\beta}{n!}\prod_{k=1}^nt_{i_k\rho_k, j_k}\langle\prod_{k=1}^n \overline{\phi}_{i_k\rho_k}\bar{\psi}_k^{j_k} \rangle_{g,n,\beta}^{\clG}.& 
\end{eqnarray}
The descendant potential $\mathcal{F}^g_{\clG}$ is a formal power series in variables $t_{i\rho, j}, 1\leq i\leq \text{rank}H^*(X,\mathbb{C}), \rho\in \widehat{G}, j\geq 0$ with coefficients in the Novikov ring $\mathbb{C}[[\overline{NE}(X)]]$, where $\overline{NE}(X)$ is the Mori cone of the coarse moduli space of $\clG$. Here $Q^\beta$ are formal variables labeled by classes $\beta\in \overline{NE}(X)$. See e.g. \cite{ts} for more discussion on descendant potentials for orbifold Gromov-Witten theory. 

Similarly the genus $g$ descendant potential of $X$ is defined to be 
\begin{equation}\label{genus0_potential_of_X}
\mathcal{F}^g_X(\{t_{i,j}\}_{1\leq i\leq \text{rank}H^*(X,\mathbb{C}), j\geq 0};Q):=\sum_{\overset{n\geq 0, \beta\in H_2(X, \mathbb{Z})}{i_1,...,i_n; j_1,...,j_n}}\frac{Q^\beta}{n!}\prod_{k=1}^nt_{i_k, j_k}\langle \prod_{k=1}^n \overline{\phi}_{i_k}\bar{\psi}_k^{j_k}\rangle_{0,n,\beta}^X.
\end{equation}
The descendant potential $\mathcal{F}^g_X$ is a formal power series in variables $t_{i,j}, 1\leq i\leq \text{rank}H^*(X, \mathbb{C}), j\geq 0$ with coefficients in $\mathbb{C}[[\overline{NE}(X)]]$  and $Q^\beta$ is (again) a formal variable.
Theorem \ref{GW-inv2} may be restated as follows.

\begin{thm}\label{decomp_any_genus}
$$\mathcal{F}^g_{\clG}(\{t_{i\rho, j}\}_{1\leq i\leq \text{rank}H^*(X,\mathbb{C}), \rho\in \widehat{G}, j\geq 0}; Q)=r^{2g-2}\sum_{[\rho]\in \widehat{G}}\mathcal{F}^g_X(\{t_{i\rho,j}\}_{1\leq i\leq \text{rank}H^*(X,\mathbb{C}), j\geq 0};Q_\rho),$$
where $Q_\rho$ is defined by the following rule: $$Q_\rho^\beta:=Q^\beta \chi_\rho\left(\exp\left(\frac{-2\pi \sqrt{-1} (\alpha\cap\beta)}{r}\right)\right),$$
and $\chi_\rho$ is the character associated to the representation $\rho$.
\end{thm}

Theorem \ref{decomp_any_genus} confirms the decomposition conjecture for genus $g$ Gromov-Witten theory of $\clG$.

\appendix
\section{A few useful facts about degree of morphisms}\label{morph-deg-appendix}
We want to characterize the degree of proper and quasi-finite morphisms of algebraic stacks in terms of
number of points in the generic geometric fibers. 
Let  $f:\clX\to\clY$ be a separated quasi-finite dominant morphism between Deligne-Mumford stacks with  $\clY$ integral.
Then according to  \cite{Vist89} its degree is given by
\begin{eqnarray}\label{deg-formula}
\mbox{deg}(\clX/\clY)=\frac{\delta(\clY)}{\delta(\clX)}\ [K(\clX):K(\clY)],
\end{eqnarray}
where  $\delta{\clX}:=\mbox{deg}(I\clX/\clX)$,  $\delta{\clY}:=\mbox{deg}(I\clY/\clY)$.
If $f: \clX\to \clY$ is representable, then the degree of the field extension  $[K(\clX):K(\clY)]$ is equal to the number of points counted with multiplicity 
in the generic geometric fiber. A similar characterization can be given in the general case for a representable  morphism 
$\overline{f}: \overline{\clX}\to \clY$
canonically determined by $f$ as described below. 
Let $\clX$ and $\clY$ be algebraic stacks. Let $f:\clX\to\clY$ be 
a  morphism  locally of finite presentation such that the relative inertia $$\mbox{Ker}(I(\clX)\to f^* I(\clY))\simeq \clX\times_{\clX\times_\clY \clX}\clX$$ is finite. Then in \cite{AOV08} it is proved that there exists a factorization of $f$
\begin{eqnarray}
\clX\stackrel{\pi}{\to} \overline{\clX}  \stackrel{\overline{f}}{\to} \clY,
\end{eqnarray}
called the {\em relative coarse moduli space}, such that 
\begin{enumerate}
\item $\overline{f}$ is representable;
 \item $\pi$ is proper  and quasi-finite;
\item $\pi_*\clO_\clX =\clO_{\overline{\clX}}$;
\item  for any $\clX\stackrel{\pi'}{\to} \clX'  \stackrel{f'}{\to} \clY$ with $f'$ representable there is a unique morphism $h:\clX'\to \overline{\clX} $ such that $\pi'=h\circ \pi$ and $\overline{f}= f'\circ h$.
\end{enumerate}
Moreover the formation of the relative coarse moduli space commutes with flat  base change and, if $f$ is tame, with arbitrary base change.

We observe that $\pi_* K(\clX)= K(\overline{\clX} )$ if $\clX$ is reduced. Indeed $K(\clX)$ is the ring of morphisms $\clX\to \mathbb{A}^1_\bbC$ defined over  an open  dense substack. For any such a morphism
there is a factorization $\clX\to \overline{\clX} \to X\to \mathbb{A}^1_\bbC$, where $\clX$ is the absolute   coarse moduli space of $X$.  Since we work over the complex numbers  and $\overline{f}:\overline{\clX}\to\clY$ is a representable morphism $[K(\overline{\clX}):K(\clY)]$ can be computed as  the number   of points (weighted by the geometric multiplicity) in a generic geometric fiber.

\bibliographystyle{alpha}

\end{document}